\documentclass[final]{siamltex}
\usepackage{float}
\usepackage{caption}
\usepackage{subfigure}
\usepackage{gensymb} 

\usepackage{amsfonts,amssymb}
\usepackage{amsmath}
\usepackage{color}
\definecolor{myblue}{rgb}{0.2, 0.5, 0.8}

\usepackage[pagewise]{lineno}
\usepackage{color}
\definecolor{darkgreen}{rgb}{0.1,0.5,0.1}
\definecolor{darkblue}{rgb}{0.2,0.2,1.0}
\usepackage[colorlinks=true,linkcolor=darkblue,citecolor=darkblue,
            filecolor=darkblue,urlcolor=darkgreen]{hyperref}

\usepackage[colorinlistoftodos]{todonotes}

\newcommand{\p}{\partial}
\newcommand{\fr}{\frac}
\def\A{{\cal A}}
\def\B{{\cal B}}
\def\W{{\cal W}}
\def\T{{\cal T}}
\def\R{{\cal R}}

\newenvironment{mat}{\left[ \begin{array}{ccccccccccccccc}}{\end{array}\right]}
\newenvironment{rmat}{\left[ \begin{array}{rrrrrrrrrrrrr}}{\end{array}\right]}
\def\bcm{\begin{mat}}
\def\ecm{\end{mat}}
\def\brm{\begin{rmat}}
\def\erm{\end{rmat}}
\newenvironment{pvect}{\left( \begin{array}{c}}{\end{array}\right)}
\def\bpvect{\begin{pvect}}
\def\epvect{\end{pvect}}

\newenvironment{choice}{\left\{ \begin{array}{ll}}{\end{array}\right.}

\def\choose{\begin{choice}}
\def\endch{\end{choice}}

\def\bsplit{\begin{split}}
\def\esplit{\end{split}}

\newcommand{\Fig}[1]{Figure~\ref{fig:#1}}

\newcommand{\ignore}[1]{}
\newcommand{\Comment}[1]{}

\newcommand{\eq}{\begin{equation}}
\newcommand{\en}{\end{equation}}
\newcommand{\eqm}{\begin{eqnarray}}
\newcommand{\enm}{\end{eqnarray}}
\newcommand{\eqmno}{\begin{eqnarray*}}
\newcommand{\enmno}{\end{eqnarray*}}
\newcommand{\eqml}[1]{\eql{#1}\begin{array}{rcl}}
\newcommand{\enml}{\end{array}\en}
\newcommand{\eql}{\begin{equation}\label}
\newcommand{\eqsub}[1]{\begin{subequations}\label{#1}\eqm }
\newcommand{\ensub}{\enm\end{subequations}}

\def\bc{\begin{center}}
\def\ec{\end{center}}
\def\bi{\begin{itemize}}
\def\ei{\end{itemize}}
\def\be{\begin{enumerate}}
\def\ee{\end{enumerate}}

\def\reals{{{\rm l} \kern -.15em {\rm R} }}

\def\qquad{\quad\quad}

\def\sf{{\cal F}}

\title{\textbf{Numerical methods for interface coupling of compressible and almost incompressible media.}
\thanks{This work was supported in part by NSF grant DMS-1216732 (RJL, MdR) and
National Council of Science and Technology of Mexico [CONACyT] (MdR).}}

\author{M.\ J.\ Del Razo\thanks{Department of Applied Mathematics, 
        University of Washington, Seattle, WA 98195-3925 ({\tt maojrs@uw.edu, rjl@uw.edu}).}
        \and R.\ J.\ LeVeque \footnotemark[2]}
\begin{document}
\maketitle
\begin{abstract}
Many experiments in biomedical applications and other disciplines
use a shock tube. These experiments often involve
placing an experimental sample within a fluid-filled container, which is 
then placed inside the shock tube. The shock tube produces an initial
shock that propagates through gas before hitting the container with the 
sample. In order to gain insight into the shock dynamics that is hard to 
obtain by experimental means, computational simulations of the 
shock wave passing from gas into a thin elastic solid and into a nearly incompressible 
fluid are developed. It is shown that if the solid interface is very thin, it 
can be neglected, simplifying the model. The model uses Euler equations for compressible 
fluids coupled with a Tammann equation of state (EOS) to model both 
compressible gas and almost incompressible materials. A three-dimensional (2D axisymmetric) 
model of these equations is solved using high-resolution shock-capturing methods, with
newly developed Riemann solvers and limiters. The methods are extended to work on a 
mapped grid to allow more complicated interface geometry, and they are adapted to work 
with adaptive mesh refinement (AMR) for higher resolution and faster computations. 
The Clawpack software is used to implement the method. These methods were initially 
inspired by shock tube experiments to study the injury mechanisms of traumatic brain 
injury (TBI).
\end{abstract}



\pagestyle{myheadings} \markboth{Interface coupling of 
compressible and almost incompressible fluids.}{M. J. Del Razo and R. J. LeVeque}
\thispagestyle{plain}

\section{Introduction}
A recent collaboration with experimentalists studying traumatic brain injury (TBI) 
at the Seattle Veterans Administration (VA) Hospital 
brought to our attention the need for very specific numerical methods \cite{delrazo2015_01}. 
Many experiments performed by the TBI community, as well as in other biomedical disciplines, employ a shock-tube, where 
they introduce samples to be studied after being exposed to a shock wave. These samples can vary from 
transwells filled with aqueous solution and cell culture to live mice \cite{coisne2005mouse,delrazo2015_01,
goeller2012investigation, goldstein2012chronic,gupta2013mathematical, huber2013blast, hue2013blood, rubovitch2011mouse}.
Within the shock-tube, the shock wave travels through highly compressible gas before hitting the 
sample, typically a nearly incompressible material with a fixed location in space. 
The physical effects of the shock wave hitting the sample are not usually evident from experimental 
data nor easy to obtain through experimental techniques. 

The methods presented in this paper were motivated by this application, although they may be useful in other
contexts as well. In order to successfully model 
the shock wave/sample interaction, we develop numerical methods that can couple the shock wave dynamics in compressible gas 
with almost incompressible materials, like plastic, water or even bone and brain. 
Some of these methods have already been employed in our recent collaboration \cite{delrazo2015_01}.
In the present paper, we give a detailed explanation of the numerical methods and their implementation; 
we extend them to more complicated interface geometries, enhance their stability
in highly refined grids, improve their resolution and 
efficiency using adaptive mesh refinement (AMR), and further study their convergence. This work will refer to 
\cite{delrazo2015_01} in the sections where it is relevant.
Although our simulations can only model idealized scenarios, they can help provide detailed insight into the behavior of 
the shock wave interaction with interfaces. For instance,
in our previous work \cite{delrazo2014,delrazo2015_01}, we obtain the dynamics of a shock wave impacting an interface that
models a specific TBI experiment. It also strongly suggested cavitation as a possible damage mechanism, an issue that
has been a subject of extense study among the TBI community \cite{delrazo2015_01,goeller2012investigation,
liu2004isentropic,moore2009computational,moss2009skull,nyein2010silico, panzer2012development,przekwas2009mathematical,ziejewski2007selected}.

Although there is an extensive body of work on computational fluid dynamics with interfaces that is
relevant and might be applicable to this type of problems, such as \cite{chertock2008interface,fagnan2008high,
gupta2013mathematical,liu2003ghost,pelanti2014mixture,saurel2001multiphase,saurel2009simple,udaykumar2002interface,wang2008adaptive} among others,
the novel methods presented here are tailored to specifically model a set of experiments performed with a shock tube. 

The methods presented here are based on finite volume methods for hyperbolic problems in their wave propagation form 
\cite{randysrbook} and implemented into Clawpack 5.2.2 \cite{clawpack}. The key ingredient in these methods 
is the Riemann solver, which must be specifically designed to deal with 
highly nonlinear waves interacting with interfaces between materials having very different properties.

In Section \ref{sec:TBInumimp}, we present the Euler Equations coupled
with the Tammann Equation of State (EOS), which can be used to model the various materials involved.
In Section \ref{sec:nummethods}, we develop the one-dimensional numerical method in detail, also discussing
its implementation in Clawpack \cite{clawpack}. The exact solution in this section is also relevant for the verification
studies in Section \ref{sec:Verif}. In Section \ref{sec:2dnummethods}, we extend the numerical 
methods to apply them to the three-dimensional (2D axisymmetric) model and on mapped grids,
allowing more complicated interface geometries. Furthermore, 
the code was designed to employ Clawpack's adaptive mesh refinement (AMR) \cite{berger1998adaptive} to improve 
efficiency and resolution. We also discuss the inclusion of transmission-based limiters to
reduce numerical oscillations in heavily refined grids produced at the interface corner. 
In Section \ref{sec:Verif}, we summarize a verification study for the one-dimensional case \cite{delrazo2015_01} and perform 
a convergence analysis for the two-dimensional method. We also show
that modifying the original minmod limiter can further reduce the numerical oscillations. The last Section discusses 
and summarizes some of the results and utility of these methods. Appendix \ref{sec:TBI:compexp} explores the question of whether a thin plastic interface between gas and liquid
can be ignored altogether in numerical studies of shock tube 
experiments \cite{delrazo2014}.
Analysis based on the nonlinear case using the numerical methods 
from Section \ref{sec:nummethods} suggests it can be ignored.
As verification, an analysis based on exact solutions to the linear acoustics equations confirms the result.

The numerical methods and implementation details are explained in this paper; the code is available in GitHub 
with a BSD license \cite{TBI-zenodo2}. 


\section{The model}
\label{sec:TBInumimp}

We use the nonlinear compressible Euler equations for compressible inviscid flow, which allow
accurate modeling of shock wave formation and propagation.  
These equations model the conservation of mass, momentum, and energy and provide a direct connection
to temperature, which may be important for some biomedical experiments.
In this type of experiment, we are not concerned with large-scale movement 
of the fluid, so viscosity does not play an important role; therefore, employing the inviscid equations is appropriate.
In order to model different materials, we use different parameters in the equations of state (EOS) for each material, so 
we can model the different materials with the same equations. The equations are solved using the methods 
explained in Section \ref{sec:nummethods}.

An additional advantage of experiments performed in a shock tube is that they often exhibit
cylindrical symmetry along the axis that goes through the center 
of the shock tube. This simplifies the three-dimensional equations into two-dimensional axisymmetric 
Euler equations, which in cylindrical coordinates ($r$,$\theta$,$z$) take the form

\begin{gather}
\label{eq:Eulercyl}
\begin{gathered}
  \fr{\p}{\p t}
    \left[\begin{array}{c} \rho \\   \rho u_r \\    \rho u_z \\    E   \end{array} \right] 
+ \fr{\p}{\p r} 
    \left[\begin{array}{c} \rho u_r \\   \rho u_r^2 + p \\    \rho u_r u_z \\     u_r(E+p)   \end{array} \right] + 
 \fr{\p}{\p z}
    \left[\begin{array}{c} \rho u_z \\   \rho u_r u_z \\    \rho u_z^2 + p \\     u_z(E+p)   \end{array} \right] 
=   \left[\begin{array}{c} -(\rho u_r)/r \\   -(\rho u_r^2)/r  \\    -(\rho u_r u_z)/r \\     -u_r(E+p)/r   \end{array} \right] ,
\end{gathered}
\end{gather}
where $\rho$ is the density; $u_r$ and $u_z$ denote the velocities in the radial and axial direction, $r$ and $z$ respectively; 
$E$ is the total energy and $p$ is the pressure. These equations have the same form as the two-dimensional Euler equations
with the addition of geometrical source terms (the right hand side), and are discussed further in Section
\ref{sec:2dnummethods}.  

\subsection{Tammann equations of state}
\label{sec:EOS}
The system of equations (\ref{eq:Eulercyl}) is closed with the addition of an EOS. It is usually given as a relation 
between pressure, density and specific internal energy, i.e. $p=p(\rho,e)$.  The most well known EOS is the 
one for an ideal gas $p = (\gamma - 1) \rho e$, where $\gamma$ is the ratio of heat capacities.
While this EOS is very good for describing the behavior of most gases, it is not appropriate 
for modeling nearly incompressible materials like water or some elastic solids.

Several alternatives exist; in this work, we will use the stiffened gas EOS, also
known as the Tammann EOS. This equation of state is very useful to model a wide range of fluids even in the 
presence of strong shock waves \cite{Fagnan2012}. The Tammann EOS is given by
\begin{align}
  p = (\gamma - 1) \rho e -\gamma p_{\infty},
  \label{eq:SGEOS}
\end{align}
where $\gamma$ and $p_{\infty}$ can be determined experimentally for different materials. 
The internal energy $e$ is related to the total energy $E$ by $E=\rho e + \frac{1}{2}\rho u \cdot u$. 
The Tammann EOS and the ideal gas EOS are the same except for the extra term $-\gamma p_{\infty}$, where
$\gamma,~p_\infty > 0$. For fluids with $p_\infty \gg p_{atm}$ (atmospheric pressure), the relative change in density, when changing 
the pressure, is very small. Consequently, the Tammann EOS is a good approximation for nearly 
incompressible fluids and can also be used to model acoustic waves in some elastic solids, like 
plastic. For sufficiently weak shocks the Tammann EOS can be further simplified to the Tait 
EOS, see \cite{Fagnan2012}, but for greater generality we use the Tammann EOS. Table \ref{tab:param} shows the Tammann EOS parameters for the
materials used in the simulations here presented.

\begin{table}[H]
  \centering
  \begin{tabular}{ l || c  c }
     Material                & $\gamma$     & $p_\infty (GPa)$ \\ \hline 
     Air (Ideal gas EOS)     & 1.4          & 0.0             \\ 
     Plastic (polystyrene)   & 1.1          & 4.79            \\ 
     Water                   & 7.15         & 0.3             \\
  \end{tabular}
  \caption{Parameters for the Tammann EOS to model the different materials. Note 
  the $p_\infty$ values for plastic and water are in GPa and are several orders of 
  magnitude above the atmospheric pressure.
  The parameters for air and water were taken from \cite{Fagnan2012}. As polystyrene
  is a solid, $\gamma$ was chosen to be close to 1, and $p_\infty$ was adjusted to yield 
  the right speed of sound in polystyrene \cite{mark2013encyclopedia}.}
  \label{tab:param}
\end{table}

\section{Numerical methods}
\label{sec:nummethods}

The Euler equations are a nonlinear hyperbolic system of conservation laws, so they can be efficiently 
solved with high-resolution shock-capturing finite volume methods (FVM). This is done by using the wave propagation algorithms 
described in \cite{randysrbook} and implemented in Clawpack \cite{clawpack}. The fundamental problem 
to solve at each cell interface of our computation is the well known Riemann problem. A general 
one-dimensional Riemann problem for a system of conservation laws like Euler equations can be stated as
\begin{align}
  q& + f(q)_x = 0, \label{eq:rphll}\\
  q&(x,0) = \left\{ 
  \begin{array}{l l}
   q_l & \  \text{if $x<0$} \\
   q_r & \  \text{if $x>0$}, 
  \end{array} \right. \nonumber
\end{align}

The Euler equatons in this work are solved by implementing a hybrid Riemann HLLC-exact type approximate solver for one-dimensional 
Euler equations with interfaces. This solver couples an HLLC approximate Riemann solver to an exact 
Riemann solver for the Tammann EOS and an Eulerian-Lagrangian description coupling at the interface. 
As the interfaces are represented by contact discontinuities, the HLLC solver is ideal to deal accurately 
with interface problems. Furthermore, the exact solver will serve as a reference solution to verify the numerical method. 

From the well-known solution to the Euler equations for an ideal EOS \cite{randysrbook,torosbook}, we expect our solution
will consist of two acoustic waves, the 1-wave and 3-wave (rarefactions or shocks), and a contact discontinuity, the 2-wave 
between them. The $n$-wave refers to the wave corresponding to the $n$-characteristic 
field (see \cite{randysrbook}).  This will separate our system in four states, $q_l,q_{*l},q_{*r},q_r$. The left 
state $q_l$ will be connected to the state $q_{*l}$ by a 1-wave, a shock wave or a rarefaction. The state $q_{*l}$ 
and $q_{*r}$ will be connected by a 2-wave, the contact discontinuity with equal pressure $p_*$ and velocity $u_*$ 
but different density on both sides. The states $q_{*r}$ and $q_r$ are connected by a 3-wave, which is a shock wave 
or a rarefaction.

The method can be extended to two dimensions by using dimensional splitting or transverse solvers. The geometrical source terms 
can be resolved using a splitting method \cite{randysrbook,randysbbook}. In the next paragraphs, we give an 
overview of the modified HLLC Riemann solver and the exact Riemann solver for the Tammann EOS with discontinuous 
parameters.

\subsection{A modified HLLC solver}

The HLLC (Harten-Lax-van Leer-Contact) solver is an approximate Riemann solver for Eq. (\ref{eq:rphll}). 
The main idea of the HLLC solver is, given the left and right going wave speeds $S_l$ and $S_r$ by some algorithm or approximation, 
assume a wave configuration of three waves separating four constant states. The Riemann solution to the 
one-dimensional Euler equations consists of three waves, two acoustic waves with a contact discontinuity 
in between. The approximate solution for this method will be of the form
\begin{align*}
  \tilde{q}&(x,t) = \left\{ 
  \begin{array}{l l}
   q_l & \  \text{if $\fr{x}{t} \le S_l$}\\
   q_{*l} & \  \text{if $S_l \le \fr{x}{t} \le S_*$},\\
   q_{*r} & \  \text{if $S_* \le \fr{x}{t} \le S_r$},\\
   q_r & \  \text{if $\fr{x}{t}\ge S_r$},
  \end{array} \right.
\end{align*}
where $S_*$ is the approximate wave speed of the contact discontinuity. Assuming we can obtain $S_l$ and $S_r$, we 
only need to find $q_{*l}$, $q_{*r}$ and $S_*$ to solve the problem. These quantities can be obtained
by integrating over a box in the $x,t$ plane using the Rankine-Hugoniot conditions and assuming constant pressure
and normal velocity across the contact discontinuity, see \cite{torosbook}. The desired states and contact discontinuity speed are given by
\begin{align*}
  q_{*k} = \fr{S_k q_k - f_k + p_* D}{S_l - S_*}, \ \ \ \text{with:} \ \ \ D = [0,1,S_*] , \\
  S_* = \fr{p_r - p_l + \rho_l u_l (S_l - u_l) -\rho_r u_r(S_r - u_r)}{\rho_l(S_l - u_l) - \rho_r(S_r - u_r)},
\end{align*}
where $\rho_k$, $u_k$ with $k=l,r$ are the left or right density and speed in the Euler equations \cite{torosbook}.

In order to calculate the wave speeds $S_l$ and $S_r$, we will
need to calculate the sound speed. This is where we require the EOS. A simple estimate is the one 
given by Davis \cite{torosbook} as 
\begin{align*}
  S_l = \text{min}\{u_l - c_l,u_r - c_r\} \ \ \  S_r = \text{max}\{u_l + c_l, u_r + c_r\},
\end{align*}
where $u_k$ is the normal velocity and $c_k$ is the sound speed on each side, $k=l,r$. Note that the easiest way to calculate the
speed of sound is using the EOS $p= p(\rho,e)$. It is usually given in the form,
\begin{align}
  c = \sqrt{\fr{\p p(\rho,e)}{\p \rho} \bigg|_s} = \sqrt{\fr{\p p(\rho,e)}{\p \rho} + \fr{p(\rho,e)}{\rho^2}\fr{\p p(\rho,e)}{\p e}},
  \label{eq:sos}
\end{align}
where $s$ is the entropy, and the first derivative is taken along the isentropic curve. A possible improvement is to employ 
Roe averages in wave speed estimates $S_l = \text{min}\{u_l - c_l, \tilde{u} - \tilde{c}\} \ \ \ 
S_r = \text{max}\{\tilde{u} + \tilde{c}, u_r + c_r\}$
where $\tilde{u}$ and $\tilde{c}$ are the Roe averages of the normal velocity and speed of sound respectively \cite{einfeldt1988godunov}. 
These Roe averages can be calculated using different configurations. Some might be more accurate when dealing with interfaces, as 
pointed out in \cite{Hu2009}.

\begin{figure}[H]
\centering
\includegraphics[width=0.9\textwidth]{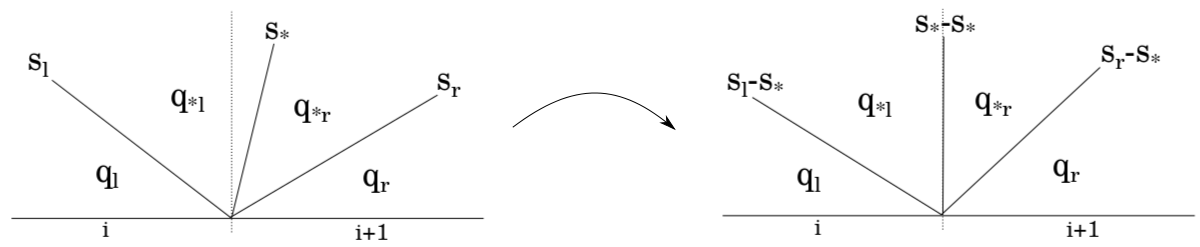}
\caption{Transformation for the HLLC Riemann solver between grid cells $i$ and $i+1$ from 
Eulerian coordinates to Lagrangian coordinates. The transformation can be employed
for other Riemann solvers too.} 
\label{fig:hllc}
\end{figure}

The HLLC solver just discussed works well for the one-dimensional Euler equations with an ideal gas EOS. 
However, we want to implement the HLLC solver with 
the Tammann EOS across an air-water or air-plastic interface. The difference between the parameters for different materials 
in the Tammann EOS are of several orders of magnitude as shown in Table \ref{tab:param}. This generates instabilities in the 
HLLC solver, more so in the multi-dimensional setting. The instability is generated because
we model the interfaces as being fixed in space; however, there is always a displacement of the contact discontinuity, i.e. 
the interface, even when the material is almost incompressible. The displacement is very small indeed, but it is big enough to 
render our numerical method unusable. In order to solve this issue, we model each material in Eulerian coordinates using the usual 
HLLC solver; if any of the cells is next to the interface, we modify our original HLLC or exact solver to work in Lagrangian 
coordinates, where the interface is actually fixed with respect to the reference frame. This is done by displacing the frame of 
reference by $S_*$, 
\begin{equation}
\label{eultolag}
\begin{aligned}
 \tilde{S_l} = S_l - S_* \ \ \
 \tilde{S_*} = 0 \ \ \
 \tilde{S_r} = S_r -  S_* 
\end{aligned}
\end{equation}
For instance, assume we are running a one-dimensional simulation of the Euler
equations, with a fixed interface modeled by a jump in the parameters of the EOS. The interface is aligned to the edge 
between cells $i$  and $i+1$, the transformed Riemann solver will be as shown in Figure \ref{fig:hllc}. This will ensure 
the contact discontinuity velocity is zero and consequently, the interface is modeled as fixed. The wave contributions will 
be the correct ones since we are just modifying the wave velocity and not the solution $q$'s. There is, of course, an error 
made at the interface when coupling the two descriptions; however, as the displacements of the interface are very small due 
to very low compressibility, this error is small, and it doesn't cause instabilities as before.


In order to provide better accuracy along the interface, we will also implement an exact Riemann solver 
for the Tammann EOS. The HLLC solver will be used to model each of the materials 
in Eulerian coordinates, and the exact solver will be used to solve the Riemann problems at the 
interface. The transformation to Lagrangian coordinates for the exact solver is equivalent to the one in (\ref{eultolag}).

\subsection{Exact Riemann solver for Tammann EOS with a jump in the parameters}
The Riemann problem (\ref{eq:rphll}) sometimes can also be solved exactly; the form of the solution will 
depend on the equations and the EOS being used. An exact solver for the Euler equations coupled with the 
Tammann EOS for constant parameters was given by Ivings \& Toro \cite{Ivings1998}. In the next paragraphs, we 
obtain the exact Riemann solver for the Euler equations coupled with the Tammann EOS with different constant 
parameters on the left and right states. This can be extended numerically to general varying parameters, by 
averaging them on each cell and using this Riemann solver to provide the solution. The solver is based on 
the one provided in \cite{Ivings1998}; however, it extends it to include a jump 
in the Tammann EOS parameters between the left and right state. 

We consider the one-dimensional Riemann problem for the Euler equations with the Tammann EOS. We want to 
solve the one-dimensional Euler equations,
\begin{align}
\label{eq:Eulercyltamman}
    \left[\begin{array}{c} \rho \\   \rho u \\    E   \end{array} \right]_t +
    \left[\begin{array}{c} \rho u \\   \rho u^2 + p \\     u(E+p)   \end{array} \right]_x
=   0,
\end{align}
where $\rho$ is density, $u$ velocity, $E$ the internal energy and $p$ the pressure and the subcripts $x,t$ 
denote partial derivatives with respect $x$ and $t$. The Tammann EOS is given by 
$p=\rho e(\gamma_k -1) - \gamma_k p_{\infty k}$ where $e$ the specific internal energy 
and $k =l,r$ determines which coefficients to use for the EOS. The initial conditions are given by the 
left and right constant states $q_l=[\rho_l, \rho_l u_l, E_l]$ and $q_r=[\rho_r, \rho_r u_r, E_r]$. Note that the 
state of the system can also be written in terms of the primitive variables $[\rho, u, p]$ by using the equation of state. 

As we mentioned before, the solution of the Euler equations will consist of the 1-wave and 3-wave 
(rarefactions or shocks), and a contact discontinuity, the 2-wave 
between them. The system will have four different solution states, $q_l,q_{*l},q_{*r},q_r$ separated by 
the three waves. In order to figure out if the 1-wave and 3-wave are rarefactions or shocks, we will need to create a function of the 
middle state pressure $p_*$ that ensures the velocity $u_*$ across the contact discontinuity is consistent. As we 
know the velocity on the left state $u_l$ should be connected by a rarefaction or shock to $u_*$, we can 
calculate $u_*=u_l + [u]_1$ with $[u]_1$the jump of the velocity across the 1-wave. In a similar manner, we also 
know the 3-wave should be a shock or rarefaction, so we can calculate $u_*= u_r - [u]_3$; therefore, we define
\begin{equation}
  \label{CDvel}
  \begin{aligned}
  \phi_l(p_*) = u_* = u_l - \mathcal{F}_l(p_*), \\
  \phi_r(p_*) = u_* = u_r + \mathcal{F}_r(p_*).
  \end{aligned}
\end{equation}
where $\mathcal{F}_{l,r}(p_*) = -[u]_{1,3}$ will change form depending if it's a shock or a rarefaction (signs were 
chosen for notation consistency). As we expect these two equations yield the same contact discontinuity 
velocity $u_*$, then 
\begin{align}
 \Phi(p_*)= \phi_r(p_*) - \phi_l(p_*) = 0.
 \label{eq:phi}
\end{align}
This nonlinear equation for $p_*$ will yield the pressure $p_*$ that provides consistency between the type of waves 
(rarefactions or shocks), their speeds and the contact discontinuity velocity $u_*$. As we mentioned before, the 
shape of $\phi_k(p_*)$ will depend on whether the states are connected by a shock wave or rarefaction. Once 
the $p_*$ has been found, the contact discontinuity velocity can be found from (\ref{CDvel}). The only remaining 
quantity to calculate from the primitive variables is the density. Furthermore, we also need the speeds of the 1-wave 
and 3-wave. Once we write the explicit equations for our system, it will be clear how to obtain these quantities. 

Before writing the equations explicitly, we should first note that having a rarefaction or shock in the 1-wave 
and 3-wave will depend on the pressure $p_*$. How can we know which one, can be answered by simple physical 
intuition. If the pressure is higher on the side toward which the wave is propagating, it will yield a rarefaction. If 
the pressure is lower, it will be a shock. In the Euler equations, this yields four possible cases for the value 
$\Phi(p_*)$ of equation (\ref{eq:phi}), just as in the solution using the ideal gas EOS \cite{Ivings1998,randysrbook}:
\begin{itemize}
 \item 1-rarefaction, 3-rarefaction: $p_*< p_l$ and $p_*<p_r$ \\
 $\Phi(p_*)= \phi_r^R(p_*) - \phi_l^R(p_*)$,
 \item 1-shock,       3-rarefaction $p_l \le p_* \le p_r$ \\
 $\Phi(p_*)= \phi_r^R(p_*) - \phi_l^S(p_*)$,
  \item 1-rarefaction,       3-shock $p_r \le p_* \le p_l$ \\
 $\Phi(p_*)= \phi_r^S(p_*) - \phi_l^R(p_*)$,
  \item 1-shock, 3-shock: $p_* > p_l$ and $p_*>p_r$ \label{4thcase}\\
 $\Phi(p_*)= \phi_r^S(p_*) - \phi_l^S(p_*)$,
\end{itemize}
where the index $S,R$ indicates if the $\phi$ was obtained by using the Rankine-Hugoniot equations to connect 
states by shocks or the Riemann invariants to connect them by rarefactions respectively.

In the next paragraphs, we derive the functions $\phi_{k}^{\mu}$
for all the four cases with $k=l,r$ and $\mu=R,S$. We show how to obtain the density and the missing wave speeds. 
In order to do so, we employ the Rankine-Hugoniot equations and the Riemann invariants. We will denote the speed 
of the 1-wave, $S_l$, the 2-wave, $S_*$, and the 3-wave $S_r$.

\subsubsection{Rankine-Hugoniot conditions for shock waves}
As we know the 1-wave and the 3-wave could each be a shock. In that case, the velocity of the wave, i.e. the 
shock, will be given by the Rankine-Hugoniot conditions. We will generalize this method for the 1-wave 
velocity $S_l$ and the 3-wave velocity $S_r$, by employing $S_k$, with $k=l,r$.

The Rankine-Hugoniot conditions are in general given by $S_k\left(q_k - q_{*k}\right) = f(q_k) - f(q_{*k})$, 
where $q$ is the vector state variable, $f(q)$ the vector state flux and $S_k$ the shock velocity. For 
the Euler equations this can be easily rewritten as \cite{Ivings1998},
\begin{gather}
 \rho_k \omega_k = \rho_{*k} \omega_*, \label{RHmass} \\
 \rho_k \omega_k^2 + p_k = \rho_{*k} \omega_*^2 + p_{*k}, \label{RHmom} \\
 \frac{1}{2} \omega_k^2 + h_k = \frac{1}{2} \omega_*^2 + h_{*k}, \label{RHene} 
\end{gather}
where $k=l,r$, $\omega_k = u_k - S_k$ , $\omega_* = u_* - S_k$ and the specific enthalpy is given 
by $h = e+(p+p_\infty)/\rho$ with $e$ the specific internal energy that relates to the internal energy 
of our original variables by $E = \rho e + \rho u^2/2$. We will use these relations to find the 
$\phi_K^S(p_*)$ of equation (\ref{eq:phi}) and the wave speeds $S_k$. \\

\noindent
\textbf{Finding $\phi_l^S(p_*)$ and $\phi_r^S(p_*)$ and $S_l$ and $S_r$:} We can start by defining 
the mass fluxes $\mathcal{Q}_k$ for $k=l,r$ as
\begin{align}
 \mathcal{Q}_l &= \rho_l \omega_l = \rho_{*l} \omega_* \label{ql0}\\
 \mathcal{Q}_r &= -\rho_r \omega_r = -\rho_{*r} \omega_*. \label{qr0}
\end{align}
As $\omega_k = u_k - S_k$, from these two equations we can obtain the wave speeds in terms of $Q_l$ and $\mathcal{Q}_r$,
\begin{equation}
\label{RHwspeeds}
\begin{aligned}
 S_l = u_l - \frac{\mathcal{Q}_l}{\rho_l}, \  \  \
 S_r = u_r + \frac{\mathcal{Q}_r}{\rho_r}.
\end{aligned}
\end{equation}
Though, we still need to find $\mathcal{Q}_l$ and $\mathcal{Q}_r$, so we substitute equation (\ref{ql0}) and (\ref{qr0}) 
into (\ref{RHmom}) to immediately obtain
\begin{align}
  \mathcal{Q}_l &= \frac{\tilde{p}_{*l} -\tilde{p}_l}{\omega_l - \omega_*} = \frac{p_* -p_l}{u_l - u_*} \label{ql1} \\
  \mathcal{Q}_r &= -\frac{\tilde{p}_{*r} -\tilde{p}_r}{\omega_r - \omega_*} = -\frac{p_* -p_r}{u_r - u_*}, \label{qr1}
\end{align}
where $\tilde{p}_\kappa = p_\kappa + p_{\infty \kappa}$ is defined to simplify future notation 
with $\kappa = l, *l, *r$ and $r$. Note that $\tilde{p}_{*l} \neq \tilde{p}_{*r}$ and 
that $\tilde{p}_{*k} -\tilde{p}_k=p_* -p_k$ since $p_*=p_{*l}=p_{*r}$ and $p_{\infty k} = p_{\infty *k}$  ($k=l,r)$.
Solving for $u_*$ we obtain the equations,
\begin{equation}
\label{RHphis}
\begin{aligned}
 u_* &= u_l - \frac{\tilde{p}_{*l} - \tilde{p}_l}{\mathcal{Q}_l} = \phi_l^S(p_*) \\
 u_* &= u_r + \frac{\tilde{p}_{*r} - \tilde{p}_r}{\mathcal{Q}_r} = \phi_r^S(p_*)
\end{aligned}
\end{equation}
Comparing to equations (\ref{CDvel}), we notice $\mathcal{F}_k(p_*) = \frac{p_* - p_k}{\mathcal{Q}_k}$. We also notice we have almost
obtained the $\phi$ functions we are looking for, though we still need to find $\mathcal{Q}_l$ and $\mathcal{Q}_r$ in terms of known 
variables.\\

\noindent
\textbf{Finding $\mathcal{Q}_l$ and $\mathcal{Q}_r$:} From equations (\ref{ql0}) and  (\ref{ql1}), we know
that
\begin{align*}
\frac{\tilde{p}_{*l} -\tilde{p}_l}{\omega_l - \omega_*} = \rho_l \omega_l.
\end{align*}
Solving for $\omega_*$, substituting the solution into (\ref{ql0}) and subtituting the $w_l$ for $\mathcal{Q}_l/\rho_l$, 
we obtain a new equation that we can solve for $\mathcal{Q}_l$ that yields,
\begin{align}
 \mathcal{Q}_k &= \sqrt{\rho_k \rho_{*k} \frac{\tilde{p}_{*k} -\tilde{p}_k}{\rho_{*k} - \rho_k}}. \label{eq:qk}
\end{align}
with $k=l,r$, since we repeated the same process for $\mathcal{Q}_r$ and obtained exactly the same equation. 
However, we still don't know $\rho_{*k}$, for this we will need our third Rankine-Hugoniot condition (\ref{RHene}).\\

\noindent
\textbf{Finding $\rho_{*l}$ and $\rho_{*r}$:} From equation (\ref{RHene}), we can obtain
\begin{align}
 h_{*k} - h_k &= \frac{1}{2}\left(w_k^2 - w_*^2\right) \nonumber\\
              &= \frac{1}{2}\left(\pm \frac{\mathcal{Q}_k^2}{\rho_k^2} \mp \frac{\mathcal{Q}_k^2}{\rho_{*k}^2}\right) \nonumber \\
              &= \frac{1}{2}\left(\frac{1}{\rho_k} + \frac{1}{\rho_{*k}}\right)(\tilde{p}_{*k}-\tilde{p}_k), \label{enthdiff}
\end{align}
where the sign above is used for $k=l$ and the one below for $k=r$, and we used equations 
(\ref{ql0}) and (\ref{qr0}) for the second line and (\ref{eq:qk}) for the third line. We can now substitute 
the specific enthalpy $h=\gamma \tilde{p}/(\rho(\gamma -1))$
in equation (\ref{enthdiff}) to obtain,
\begin{align*}
 \frac{\gamma_{*k}}{\gamma_{*k} - 1}\frac{\tilde{p}_{*k}}{\rho_{*k}} - \frac{\gamma_{k}}{\gamma_{k} - 1}\frac{\tilde{p}_k}{\rho_{k}}
 =\frac{1}{2}\left(\frac{1}{\rho_k} + \frac{1}{\rho_{*k}}\right)(\tilde{p}_{*k}-\tilde{p}_k).
\end{align*}
As the interface is the contact discontinuity, the jump in the parameters is only across the contact discontinuity, so 
$\gamma_{*k}=\gamma_k$. Now we can solve for the unknown
density,
\begin{align}
 \rho_{*k} = \rho_k \left(\frac{\frac{\tilde{p_*}}{\tilde{p_k}} + 
 \frac{\gamma_k-1}{\gamma_k+1}}{\frac{\tilde{p_*}}{\tilde{p_k}}\frac{\gamma_k-1}{\gamma_k+1}+1}\right). 
 \label{rhok}
\end{align}
Replacing this result into equation (\ref{eq:qk}), we obtain $\mathcal{Q}_k$ in terms of $p_*$ and known variables,
\begin{align}
\label{qkfinal}
 \mathcal{Q}_k = \sqrt{\rho_k \frac{\tilde{p}_{*k} + \tilde{p}_{k}\frac{\gamma_k - 1}{\gamma_k +1}}{\frac{2}{\gamma_k+1}}}.
\end{align}
With equations (\ref{RHphis}) and (\ref{qkfinal}), we can calculate the $\phi^S_{l,r}$ nonlinear functions of $p_*$
in terms of known variables. The functions $\phi^S_{l,r}$ allow us to construct equation (\ref{eq:phi}) and solve it using 
a Newton method 
or other root finder in order to obtain the value of $p_*$.
Equations (\ref{RHphis}) will then yield the contact discontinuity speed $S_*=u_*$ in terms of $p_*$. 
Further on, we can calculate $\mathcal{Q}_l$ and $\mathcal{Q}_r$ from (\ref{qkfinal}), and we can substitute in (\ref{RHwspeeds}) 
to obtain the corresponding wave speeds. However, this will only solve the 4th case of equation (\ref{eq:phi}), 1-shock and 
3-shock solution. If any of our waves happens to be a rarefaction, we will also need to calculate the $\phi^R_{l,r}$ 
functions. This will be obtained using the Riemann invariants.

\subsubsection{Riemann invariants for rarefaction waves}
Riemann invariants are variables that remain constant through simple waves such as rarefactions. The Riemann invariants across 
the 2-wave are the pressure $p_*$ and the normal velocity $u_*$. The Riemann invariants for the 1-wave and 3-wave are 
the entropy and the quantities,
\begin{align}
 u_l+\frac{2c_l}{\gamma_l-1} = u_* + \frac{2c_{*l}}{\gamma_l-1}, \label{RI1} \\
 u_r-\frac{2c_r}{\gamma_r-1} = u_* - \frac{2c_{*r}}{\gamma_r-1}, \label{RI2}
\end{align}
correspondingly. The speed of sound $c_K$ is obtained by applying equation (\ref{eq:sos}) to the Tammann EOS,
\begin{align}
 \label{sndsp}
 c_k= \sqrt{\gamma_k\frac{p_k+p_{\infty k}}{\rho_k}}.
\end{align}
As the entropy is invariant, we can use the Tammann EOS isentropic relation to obtain the density in the middle states,
\begin{align}
 \label{isentrop}
 \rho_{*k} = \rho_k \left(\frac{\tilde{p}_{*k}}{\tilde{p}_k}\right)^{1/\gamma}.
\end{align}
Solving (\ref{RI1}) and (\ref{RI2}) for $u_*$ and using equations (\ref{sndsp}) and (\ref{isentrop}), we immediately obtain
\begin{align*}
 u_* = u_l + \frac{2c_l}{\gamma_l-1} \left[
 1-\left(\frac{\tilde{p}_{*l}}{\tilde{p}_l}\right)^{\frac{\gamma_l-1}{2\gamma_l}}
 \right] = \phi_l^R(p_*), \\
 u_* = u_r - \frac{2c_r}{\gamma_r-1} \left[
 1-\left(\frac{\tilde{p}_{*r}}{\tilde{p}_r}\right)^{\frac{\gamma_r-1}{2\gamma_r}}
 \right] = \phi_r^R(p_*).
\end{align*}
As $\tilde{p}_\kappa = p_\kappa + p_{\infty \kappa}$, when we compare to equations (\ref{CDvel}) we obtain 
the $\phi_{l,r}^R$ functions. The rarefaction head velocities will be given by $u_l-c_l$ 
and $u_r+c_r$; the tail velocities will be $u_*-c_{*l}$ and $u_*+c_{*r}$. For numerical purposes, a simple 
approximate velocity is provided for $S_l$ and $S_r$ as the average between the head and tail velocity. 

In order to compute the complete structure of the rarefaction wave \cite{Ivings1998,randysrbook}, we can use the 
Riemann invariants from Eqs. \ref{RI1} and \ref{RI2}, along with Eq. \ref{sndsp} and the isentropic 
relation from Eq. \ref{isentrop}. The solution for the 1-rarefaction wave along the rays $x/t=\xi= u_{rar1} - c_{rar1}$
is then
\begin{align*}
 u_{rar1}(\xi) &= \fr{u_l(\gamma_l - 1) + 2(\xi + c_l)}{\gamma_l + 1}, \\
 \rho_{rar1} &= \rho_l \left[\fr{u_{rar1}(\xi) - \xi}{c_l}\right]^{\fr{2}{\gamma_l -1}}, \\
 p_{rar1} &= \tilde{p}_l \left[\fr{u_{rar1}(\xi) - \xi}{c_l}\right]^{\fr{2\gamma_l}{\gamma_l -1}} - p_{\infty l},
\end{align*}
and for a 3-rarefaction wave along the rays $x/t=\xi= u_{rar3} + c_{rar3}$ is,
\begin{align*}
 u_{rar3}(\xi) &= \fr{u_r(\gamma_r - 1) + 2(\xi - c_r)}{\gamma_r + 1}, \\
 \rho_{rar3} &= \rho_r \left[\fr{u_{rar3}(\xi) - \xi}{c_r}\right]^{\fr{2}{\gamma_r -1}}, \\
 p_{rar3} &= \tilde{p}_r \left[\fr{u_{rar3}(\xi) - \xi}{c_r}\right]^{\fr{2\gamma_r}{\gamma_r -1}} - p_{\infty r}.
\end{align*}

Now that we know the functions $\phi^{s,r}_{l,r}$ for the rarefactions, we can construct the function 
$\Phi(p_*)$ function from (\ref{eq:phi}) for any of the 4 possible scenarios. The value of $p_*$ will be 
found by numerically finding the roots of $\Phi(p_*)=0$. Note which case to employ to calculate $\Phi(p_*)$ 
might change in each iteration of the root finder. Once $p_*$ is found, $u_*$, $\rho_{*l}$, $\rho_{*r}$, $S_l$ 
and $S_r$ can be found using the relations we just derived depending if it's a shock or a rarefaction. As we 
know the three wave speeds $S_l$, $S_*$ and $S_r$ and the primitive variables $[\rho,u,p]$ on all the 4 
states for all the possible cases, we have the solved the Riemann problem.

\subsection{Implementation into Clawpack}
These methods are implemented into the Clawpack 5.2.2 software \cite{clawpack}. This software employs
Godunov's method \cite{godunov1959difference} with high order corrections and limiters to better handle 
discontinuities\cite{randysrbook}. In order to implement these methods into Clawpack, we
need to write Godunov's method in the wave propagation form. Consider a state vector $q(x,t)$, a one dimensional 
conservation law is given by $q_t + f(q)_x = 0$. We partition the space in cells with index $i$ and 
consider the cell average at time $t$ to 
be $Q_i^n = \int_{x_{i-1/2}}^{x_{i+1/2}} q(x,t_n)dx $. Then the Godunov method is given by,
\begin{align}
 Q_{i}^{n+1} = Q_i^n - \frac{\Delta t}{\Delta x} (
\underbrace{\A^- \Delta Q_{i+1/2} }_{\mathrm{Left Edge}} +  
\underbrace{\A^+ \Delta Q_{i-1/2}}_{\mathrm{Right Edge}} )
- \underbrace{\frac{\Delta t}{\Delta x}\left(\tilde{F}_{i+1/2} - \tilde{F}_{i-1/2}\right)}_{
  \mathrm{High Resolution}},
  \label{eq:Godunov}
\end{align}
with,
\begin{align}\tilde{F}_{i \pm 1/2} = \frac{1}{2}\sum_{p=1}^{m} |s_{i \pm 1/2}^p |
                            \left(1 - \frac{\Delta t}{\Delta x} |s_{i \pm 1/2}^p| \right) 
                            \underbrace{\widetilde{\W}_{i \pm 1/2}^p}_{\mathrm{Limiter}},  
 \label{eq:Godunovpt2}
\end{align}
where $\A^- \Delta Q_{i\pm1/2} = \sum_{p=1}^m 
(s^p_{i\pm1/2})^- \W_{i\pm1/2}^p$ and $\A^+ \Delta Q_{i\pm1/2} = \sum_{p=1}^m 
(s^p_{i\pm1/2})^+ \W_{i\pm1/2}^p$ are the left and right going fluctuations of the edge of cell $i\pm 1/2$
respectively, with $(s^p_{i\pm1/2})^\pm$ indicating only those values of $s^p_{i\pm1/2}$ with sign $\pm$, $m$ is 
the number of waves, $s^p_i\mp 1/2$ is the velocity of the $p$ characteristic of the Riemann problem at
edge $i\mp 1/2$, the wave $\W_{i\mp 1/2}^p$ corresponds to the jump across that characteristic and $\widetilde{\W}_{i \pm 1/2}^p$
is the limited version of the wave, see \cite{randysrbook} for more details. 

The numerical solution requires solving a Riemann problem on each cell edge of our partition in order to
obtain the fluctuations. The Riemann solutions presented previously have provided the characteristic 
velocities $s^p$, and we can calculate the waves $\W^p$ by calculating the jump of $q$ across the $p$ 
characteristic. This information is calculated for each cell edge and fed into Clawpack, where the method
from Eq. \ref{eq:Godunov} is implemented. Appendix \ref{sec:TBI:compexp} contains one-dimensional implementations
of these methods. In the next Sections, we will study two-dimensional
implementations.

\section{Two dimensional axisymmetric model}
\label{sec:2dnummethods}
\begin{figure}[h]
  \centering
  \includegraphics[width=0.5\textwidth]{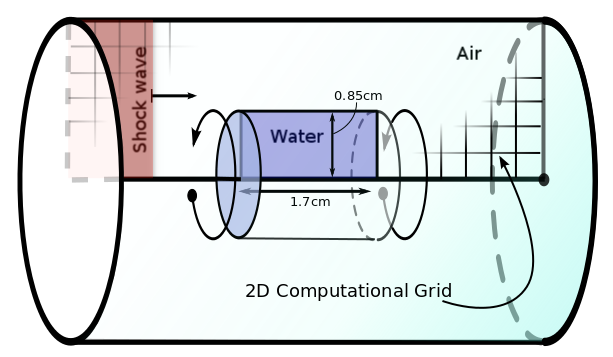} 
  \caption{The axisymmetric model is obtained by revolving the
2D computational grid. The inner square corresponds to the air-water interface. The inside part is filled with water and the outside 
part is filled with air. All the outer boundaries are modeled with non-reflecting boundary conditions. The interface location
was chosen following the source of this figure \cite{delrazo2015_01}.}
  \label{fig:2drevolve}
\end{figure}

The three dimensional Euler Equations with cylindrical symmetry can be solved as two dimensional axisymmetric Euler Equations 
with additional source terms, see Eqs. \ref{eq:Eulercyl} and \Fig{2drevolve}.
The conservation law for $q(x,y,t)$ takes the form  $q_t + f(q)_x + g(q)_y = \psi(q,x,y,t)$. In two dimensions, the 
numerical cell average is calculated as 
$Q_{i,j}^n=\fr{1}{\Delta y\Delta x}\int_{C_{i,j}} q(x,y,t_n)dxdy$, 
where $C_{i,j}$ is the cell $[x_{i-1/2}, x_{i+1/2}]\times [y_{j-1/2},y_{j+1/2}]$. The source terms can be solved using a
fractional-step method \cite{randysrbook} by alternating between $q_t + f(q)_x + g(q)_y = 0$ and $q_t = \psi(q,x,y,t)$. The latter is an
ordinary differential equation, which has an exact solution in the case of Equations \ref{eq:Eulercyl}, as shown in \cite{delrazo2015_01}. 
More complex source terms might require implementing another time stepping method like Runge-Kutta or TR-BDF2. In a 
similar manner, the simplest approach to solve the two dimensional system $q_t + f(q)_x + g(q)_y = 0$ is dimensional splitting. 
This is done again with a fractional-step method to split the two dimensional problem up into a
sequence of one-dimensional problems alternating between solving $q_t + f(q)_x = 0$ and 
$q_t + g(q)_y = 0$. For more details and different splitting algorithms see \cite{randysrbook}. 

Although dimensional splitting is simple to implement, we can obtain
second-order accuracy and less numerical smearing simultaneously by using transverse propagation algorithms from \cite{leveque1997wave}. This will require splitting the normal wave 
fluctuations $\A^\pm \Delta Q_{i\pm 1/2,j}$ at edge $i\pm1/2$ into transverse wave 
fluctuations $\B^\pm \A^+\Delta Q_{i\pm 1/2,j}$ and $\B^\pm \A^-\Delta Q_{i\pm 1/2,j}$. If the normal direction is $x$, then 
the normal fluctuations are calculated with the flux $f(q)$ and the transverse ones with the flux $g(q)$. Our
specific model will require a very special kind of transverse solvers, which have been implemented in \cite{delrazo2015_01}; a
generalized version of these solvers will be explored in detail later in this paper.

The two-dimensional axisymmetric model of Eqs. \ref{eq:Eulercyl} employing a Tammann equation of state with interfaces and
transverse solvers were implemented in a traumatic brain injury application in \cite{delrazo2015_01}. This work showed how 
the geometry of the interface can be very relevant and even produce cavitation effects. The set up in \cite{delrazo2015_01} and in this work
is essentially the one shown in \Fig{2drevolve}. A cylindrical plastic container filled with water is placed inside a shock tube. The cylindrical 
outer boundary corresponds to a cylindrical cross section of the shock tube. The results shown in the Appendix \ref{sec:TBI:compexp} and 
in \cite{delrazo2014} show the plastic interface can be neglected in the two-dimensional model.

\begin{figure}[h]
  \centering
  \includegraphics[width=0.94\textwidth]{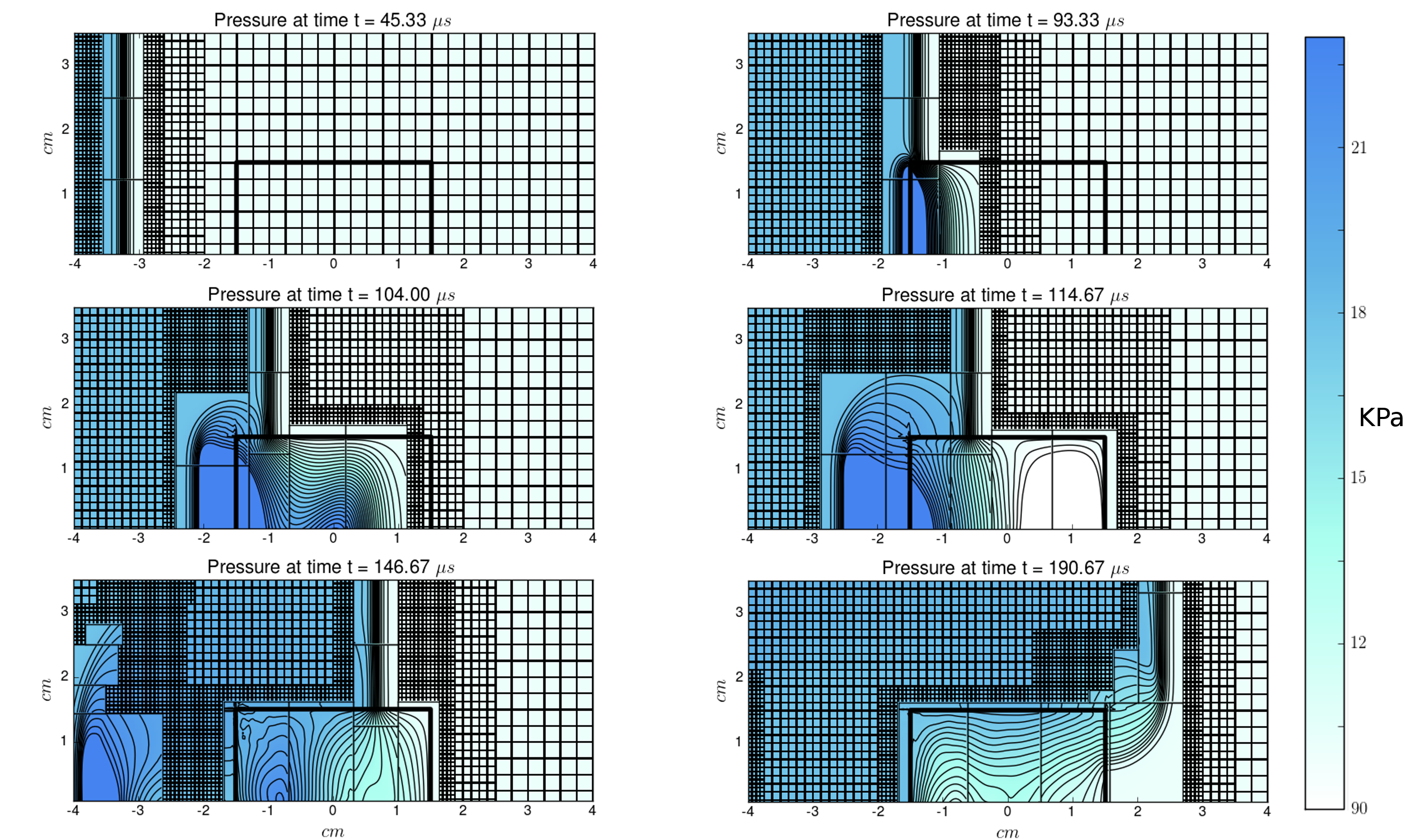} 
  \caption{Axisymmetric simulation pressure contour plots at six different times points
  $t=45.33, 93.33, 104, 114.67, 146.67, 190.67 \mu s$, using four levels of AMR. The parameters employed to model water and air
  for the Tammann EOS are the ones in Table \ref{tab:param}. The pressure amplitude is given along the color bar in KPa.
  The interface separating air and water is marked as a thick black line, and considering the axis of symmetry is the 
  $x$ axis, it models a cylindrical water interface immersed in air. The shock wave travels from the left to right. The 
  first, second and third AMR grid refinement levels are plotted explicitly while the fourth level just
  shows the refinement patches for clarity. The pressure contours are only shown in the highest refinement level.}
  \label{fig:cartesian}
\end{figure}

In this work, the model implemented in \cite{delrazo2015_01} is extended to work with AMR capabilities in 
Clawpack \cite{clawpack,berger1998adaptive}. The AMR implementation requires interpolating the value from coarser grid cells 
into the finer ones. However, when this
interpolation is done across the interface, it will cause instabilities due to the big jump in the EOS parameters
across the interface. In order to address this issue, we had to make sure that when a refinement patch intersects the interface,
the interpolation for the finest grids is performed only using grid cells corresponding to the same material. For instance,
if we need to refine a water grid cell, which is adjacent to the air interface, we will only use the values
of adjacent cells corresponding to other water grid cells to obtain the interpolated values in the refined cells. It should be noted that
the interface is always aligned to the cell edges, so there are no grid cells that contain two materials. This is also
true for the mapped grid case studied below. \Fig{cartesian} shows 
the pressure contours for six different time points
for a shock wave traveling in air and hitting a water interface fixed in space, as illustrated in \Fig{2drevolve}. The 
grid is plotted on top showing AMR in action with 4 levels of refinement. The first, second and third coarser grid 
levels are shown explicitly. The level four refinement is plotted as patches that indicate the highest refinement. 
Additionally, the code allows us to add gauges to observe the pressure as a function of time at any given point.  

In addition to the implementation of these methods in \cite{delrazo2014,delrazo2015_01}, we now show an extension of the algorithms 
for a mapped grid with adaptive mesh refinement (AMR).

\subsection{Two dimensional model in a mapped grid}
These algorithms can also be used on a mapped grid where the quadrilateral grid cells are not necessarily rectangular. 
We will first consider how to implement the normal Riemann solver in the mapped grid.
This will require a mapping from a Cartesian grid to a quadrilateral grid, which will tell us the normal at each cell edge where
we are solving the Riemann solver as well as the scaling of the edges and the scaling of the areas of the cells.
The mapped normal Riemann solver can be done using the same solver as in the Cartesian case by following these steps: 
\begin{itemize}
 \item Define a mapping;
 \item Use the normal at each mapped cell edge to rotate the velocities from the computational 
 domain into normal and transverse components in the physical domain;
 \item Solve the Riemann problem as usual with the rotated velocities and calculate the waves;
 \item Rotate the waves back into the computational domain;
 \item Use the cell edge and area scaling to modify the algorithm in Eq. \ref{eq:Godunov}, see \cite{randysrbook}.
\end{itemize}

\begin{figure}[h]
  \centering
  \includegraphics[width=0.47\textwidth]{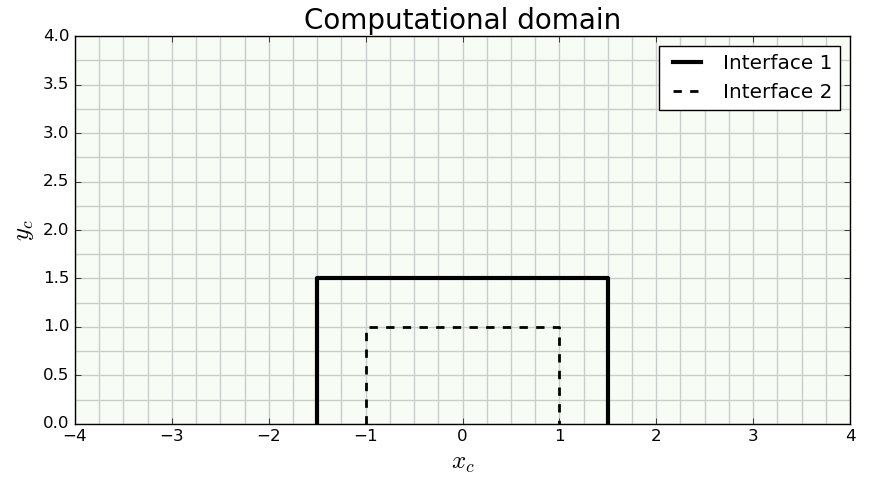} 
  \includegraphics[width=0.47\textwidth]{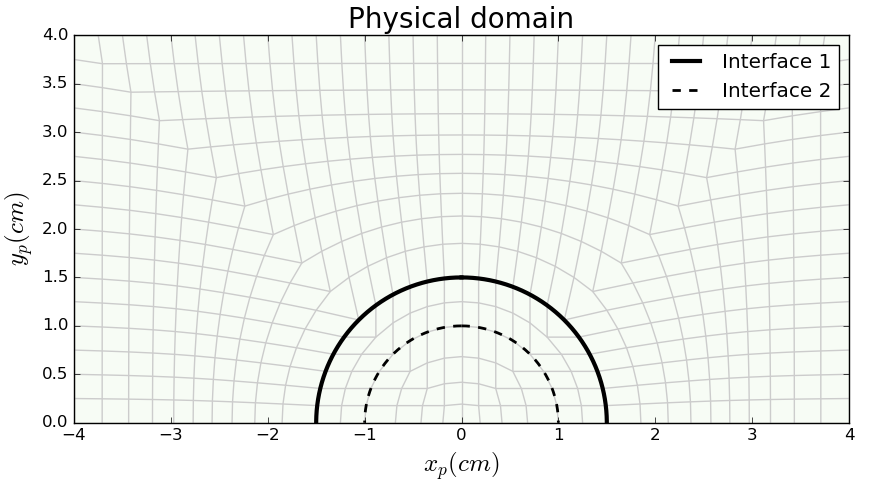} 
  \caption{Computational and physical mapped grid of a circular shell inclusion based on the mapping in 
  \cite{calhoun2008logically}. The mapping provides two possible circular
  interfaces, so considering the model is axisymmetric along the $x$ axis, it can be used to model a 
  spherical interface or a spherical thick shell interface. 
  The locations of two possible interfaces 
  are shown as thick continuous and dashed lines in both domains.}
  \label{fig:grids}
\end{figure}

The mapping of Figure \ref{fig:grids} is based on the mappings of \cite{calhoun2008logically}. Consider
a computational point $(x_c,y_c)$ on a rectangular grid
such that $x_c>0$ and $|y_c|<x_c \equiv d$. The vertical 
line segment from $(d,-d)$ to $(d,d)$ will be mapped to a circular arc with radius $R(d)$ 
that intersect the identity diagonals at $(D(d),-D(d))$ and $(D(d),D(d))$. The center of
such a circular arc is then given by $(x_0,y_0) = (D(d)-\sqrt{R(d)^2 - D(d)^2},0)$, and the 
point in the computational grid is mapped to the physical grid point $(x_p,y_p)$ by
\begin{align*}
 y_p &= y_c D(d)/d, \\
 x_p &= x_0 + \sqrt{R(d)^2 - y_p^2}.
\end{align*}
In the mapping of Figure \ref{fig:grids} we have indicated two interfaces: the inner one at radius $r_i$ (1cm) and the outer one 
at radius $r_{o}$ (1.5cm) from the origin. The size of the square domain in the computational grid where the mapping is applied is given 
by a third parameter $r_m$ (4cm), with $r_m > r_o > r_i$. The square domain is centered at the origin and the length of each side is $2r_m$. 
In order to determine the mapping, we need to choose $R(d)$ and $D(d)$ in the three regions defined by the two interfaces. One option
that works well, as shown in Figure \ref{fig:grids}, is given by
\begin{align*}
   D(d) = 
   \begin{cases} 
      r_m \fr{d}{\sqrt{2}} \\[2mm]
      r_m \fr{d}{\sqrt{2}} \\[2mm]
      \fr{r_o}{\sqrt{2}} + \fr{\left(d-\fr{r_o}{r_m}\right)\left(r_m-\fr{r_o}{\sqrt{2}}\right)}{1-\fr{r_o}{r_m}}
   \end{cases},
   \hspace{5mm}
   R(d) = 
   \begin{cases} 
      r_i & d \leq \fr{r_i}{r_m} \\[2mm]
      d r_m & \fr{r_i}{r_m} < d \leq \fr{r_o}{r_m} \\[2mm]
      r_m \left[\fr{1-\fr{r_o}{r_m}}{1-d}\right]^{\left(\fr{r_m}{r_o} + \fr{1}{2}\right)} & d > \fr{r_o}{r_m}
   \end{cases}.
\end{align*}
Note this is only for the eastern sector of the computational grid, where $x_c>0$ and $|y_c|<x_c$; the other sections 
are analogous \cite{calhoun2008logically}.

Some of the quadrilateral cells in the physical domain are nearly
triangular, with two adjacent edges nearly colinear.  In spite of this, the
wave-propagation algorithm with transverse solvers described below works quite
robustly in general as discussed further in \cite{calhoun2008logically}.
However, when there is also a large jump in material
parameters at the interface and the grids are adaptively refined there can be some
stability issues as discussed further below.

Once the mapping is defined, we proceed by rotating the normal and transverse momentum components $q^2$ and $q^3$ of 
the Euler equations in the computational grid by using the normal at the current edge of the mapped grid, $\hat{n}=(n_x,n_y)$,
 \begin{gather*}
\begin{gathered}
    \left[\begin{array}{c} q_{ph}^2 \\ q_{ph}^3 \end{array} \right] =
    \left[\begin{array}{c c} n_x & n_y \\  
                -n_y & n_x \end{array} \right]
    \left[\begin{array}{c} q^2 \\ q^3 \end{array} \right],
\end{gathered}
\end{gather*}
where $q_{ph}^2$ and $q_{ph}^3$ now point in the normal and transverse direction in the physical domain (mapped grid). Using these 
quantities, we solve the normal Riemann solver as usual to obtain the speeds and waves $s^p_{ph}$ and $\W^p_{ph}$, and we rotate
the waves back to the computational domain,
\begin{gather*}
\begin{gathered}
    \left[\begin{array}{c} \W^2 \\ \W^3 \end{array} \right] =
    \left[\begin{array}{c c} n_x & -n_y \\  
                 n_y & n_x \end{array} \right]
    \left[\begin{array}{c} \W_{ph}^2 \\ \W_{ph}^3 \end{array} \right].
\end{gathered}
\end{gather*}
Finally, we scale the speeds $s^p_{ph}$ by the edge scaling to obtain $s^p$ and employ the capacity function (cell area scaling) into 
a modified version of the algorithm in Eq. \ref{eq:Godunov} found on \cite{randysrbook}. 

The transverse solvers will also be applied on the mapped grid, but this requires more careful consideration because of our treatment of 
interfaces with huge jumps in the Tammann EOS parameters. This will be explained in detail in the next subsection. In this work,
we implemented the 2D axisymmetric model into the mapped grid of Figure \ref{fig:grids}. Although the mapping is two-dimensional and
shows half circular inclusions interfaces, the axisymmetry along the $x$ axis convert these interfaces into spherical shells. This 
mapped grid was selected because it could be used to model a skull in computational TBI 
experiments. Note the mapping allows an inner interface that could even
be used to model the thickness of a skull. The code is set up so 
arbitrary mappings with other interface geometries can be implemented. 

\begin{figure}[h]
  \centering
  \includegraphics[width=0.94\textwidth]{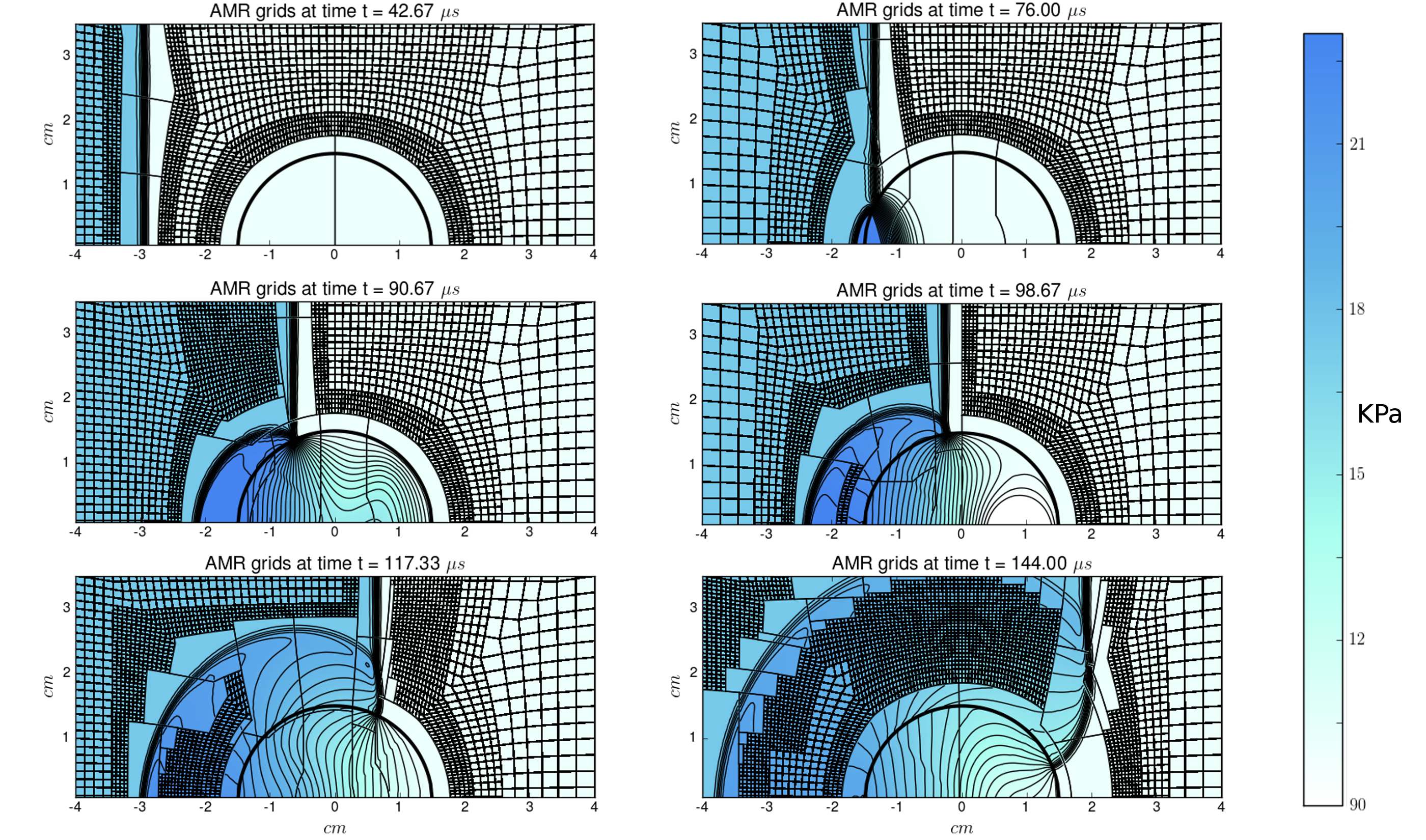} 
  \caption{Pressure contour plots of axisymmetric simulation on a mapped grid with a circular inclusion at six different times points
  $t=42.67, 76, 90.67, 98.67, 117.33, 144$ $\mu s$, using four levels of AMR. The plot is analogous to that of Figure \ref{fig:cartesian}; however,
  in this figure the interface separating air and water is circular, which models a spherical water interface. Also note, the region
  around the interface is refined from the beginning to avoid instabilities when using AMR around corners in the mapped grid.}
  \label{fig:maphighP}
\end{figure}

In Figure \ref{fig:maphighP}, we show a sample simulation of the pressure contours for the mapped grid 
at six different points in time. It only employs one interface along the outer circular inclusion shown in the grid of Figure \ref{fig:grids}.
Once again the outer part of the circular inclusion is modeled as air and the inner material as water using the same set 
of parameters as Figure \ref{fig:cartesian}. This figure also 
shows AMR in action with 4 levels of refinement, and it is also possible to add gauges to observe the pressure as a function of time at 
any given point in the grid. AMR does not need many additional considerations in terms of the mapped grid since it works on the computational
domain, which is still Cartesian. However, it is worth mentioning that the region around the interface is refined to the highest level 
from the beginning of the simulation. This is to avoid instabilities caused by employing AMR along an interface with 
huge jumps in the parameters while using a mapped grid with almost triangular grid cells (see Figure \ref{fig:grids}). If any of the 
conditions is relaxed, i.e. we use a smaller jump in the parameters or use a less severe mapped grid as in Figure
\ref{fig:cartesian}, this initial refinement along the interface is no longer required to avoid instabilities.

\subsection{Transverse Riemann solver in a mapped grid}
A transverse solver for a Cartesian grid was implemented in \cite{delrazo2015_01}. In this section, we show the extension
of this transverse Riemann solver for a mapped grid. This solver takes the results of a normal Riemann solver and splits it
into components moving in the transverse direction. As mentioned in \cite{delrazo2015_01},  
a special transverse solver needs to be developed due to instabilities at the interface. This is
based on the solver for acoustics in a heterogeneous media that is described in Section
21.5 of \cite{randysrbook}.

We recall the basic idea of a transverse solver for a constant 
coefficient linear hyperbolic system of equations 
$q_t + A q_x +B q_y =0$, the jump in normal flux between adjacent cells,
$\A\Delta Q_{i-1/2} = \A(Q_{i,j} - Q_{i-1,j})$, is split via the normal Riemann
solver into left-going and right-going ``fluctuations'' $\A^-\Delta Q_{i-1/2}$ and $\A^+\Delta Q_{i-1/2}$.  
Each fluctuation $\A^+\Delta Q_{i-1/2}$, is then
further split into down-going and up-going components $\B^-\A^+\Delta
Q_{i-1/2}$ and $\B^+\A^+\Delta Q_{i-1/2}$, based on the matrices $B^+$ and $B^-$.  

In the case of variable coefficients or nonlinear problems, 
the general notation $\B^-\A^+\Delta Q_{i-1/2}$ and $\B^+\A^+\Delta
Q_{i-1/2}$ is used for these two vectors.  For variable coefficient
acoustics, as described in \cite{randysrbook}, the up-going fluctuation from
the transverse splitting is based on eigenvectors of $B_{ij}$ and $B_{i,j+1}$,
while the down-going fluctuation is based on eigenvectors of $B_{ij}$ and
$B_{i,j-1}$.

At the interface with an almost incompressible liquid, it is difficult to figure out an accurate and stable implementation of the 
transverse Riemann problem. This is because Euler equations, with a big jump in the parameters at the interface, are extremely sensitive to 
instabilities. Our first approach was to expand the normal wave as a function of linearized eigenvectors 
corresponding to the transverse grid cells \cite{randysrbook} of the Euler equations. However, this approach resulted in instabilities
at the interface. In order to work around this issue, we will follow the same approach as \cite{delrazo2015_01} and derive an 
approximate transverse Riemann solver based on acoustic equations, which capture the acoustic waves while avoiding instabilities.

In this interface, we will mostly be concerned with the two acoustic waves.
In order to derive it, let $\hat{n}=(n_x,n_y)$ be the transverse unitary normal vector and 
linearize the acoustic equations around $\rho_0,
u_0$ and $v_0$, with $u_0$ and $v_0$ the velocity in the $x$ and $y$ direction 
respectively \cite{randysrbook}.  In terms of the density and momentum,

 \begin{gather}
\label{eq:linacoust}
\begin{gathered}
    \left[\begin{array}{c} \rho \\   \rho u \\ \rho v \end{array} \right]_t +
    \underbrace{\left[\begin{array}{c c c } 0 & n_x & n_y \\  
                        n_x c^2 & 0 & 0 \\
                        n_y c^2 & 0 & 0 \end{array} \right]}_{ \tilde{B}(Q)}
    \left[\begin{array}{c} \rho \\ \rho u \\ \rho v \end{array} \right]_{\hat{n}}
= 0,
\end{gathered}
\end{gather}
where the derivative is taken in the normal direction $\hat{n}$, $c$ is the sound speed 
and $\tilde{B}(Q)$ can be understood as a lower dimensional approximation of the transverse Jacobian 
$g'(Q_0)$ for the Euler equations. Note we assumed $u_0=0$, which is equivalent to move into a Lagrangian frame of reference. 

As we might have different materials and sound speeds in the cell above or below, we calculate the eigenvectors
and evaluate them according to their location. The matrix of eigenvectors is

 \begin{gather*}
\begin{gathered}
R=
    \left[\begin{array}{c c c} 1 & 1 & 0 \\  
                             n_{x}c_U & -n_{x} c_D & -n_y \\
                             n_{y}c_U & -n_{y} c_D & n_x 
          \end{array} \right],
\end{gathered}
\end{gather*}
where the sound speeds $c_U$ and $-c_D$ are the eigenvalues corresponding to the first two column eigenvectors, $v_u$ and $v_d$.
The eigenvalue for the third one $v_0$ is 0. The subindex $u$ and $d$ refer to cells $(i,j+1)$ and $(i,j)$
when computing $\B^+\A^+\Delta Q_{i-1/2,j}$ and to cells $(i,j)$ and $(i,j-1)$ when computing $\B^-\A^+\Delta Q_{i-1/2,j}$.

The up-going and down-going fluctuations for $\A^+\Delta Q_{i-1/2,j}$ are obtained by expanding the fluctuation 
in terms of these eigenvectors or waves, $\A^+\Delta Q_{i-1/2,j} = \alpha_U v_U + \alpha_D v_D + \alpha_0 v_0$,
so we need to solve $R\alpha=\A^+\Delta Q_{i-1/2,j}$, which yields

\begin{align*}
 \alpha_U = \frac{1}{c_U + c_D}\left( c_D \A^+_{1} + n_x \A^+_{2} + n_y \A^+_{3}  \right), \\
 \alpha_D = \frac{1}{c_U + c_D}\left( c_U \A^+_{1} - n_x \A^+_{2} - n_y \A^+_{3}  \right),
\end{align*}
and $\alpha_0$ is not relevant since it corresponds to the zero eigenvalue.
Note that the required fluctuation $\A^+\Delta Q_{i-1/2,j}$ for 
the Euler equations is a four-dimensional vector with fluctuations in density, normal momentum, transverse momentum, and 
energy. As we are only interested in the acoustic waves, we will assume the fluctuations in energy are 
negligible, so we define the acoustic part of the fluctuation as the first second and third entry of the 4 dimensional vector, 
i.e. $\A^+_{ac}\Delta Q_{i-1/2,j} = [\A^+_{1},\A^+_{2},\A^+_{3}]$. 

\begin{figure}[t]
\centering
\includegraphics[width=0.38\paperwidth]{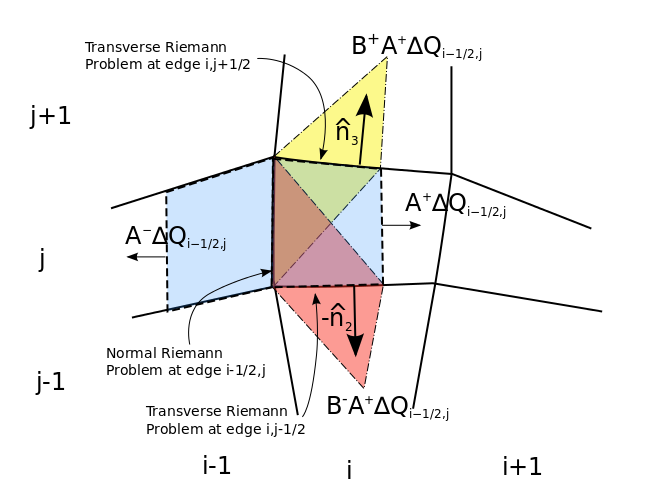}
\caption{Transverse solvers diagram in physical grid cells after applying the mapping. 
The left-going and right going fluctuations of 
the normal Riemann problem at the edge between grid cells $(i-1,j)$ and $(i,j)$ is shown. The right-going 
fluctuation $\A^+\Delta Q_{i-1/2,j}$ is decomposed into the up-going fluctuation $\B^+\A^+\Delta Q_{i-1/2,j}$
and the down-going fluctuation $\B^-\A^+\Delta Q_{i-1/2,j}$ by employing transverse Riemann solvers in the 
computational grid. This is an extension of the transverse solvers implemented in \cite{delrazo2015_01} into mapped 
grids.} 
\label{fig:trans}
\end{figure}

The up-going and down-going acoustic fluctuations are given by the velocity times the waves,
\begin{align*}
 \B^+_{ac}\A^+\Delta Q_{i-1/2,j} = c_U \alpha_U v_U, \\
 \B^-_{ac}\A^+\Delta Q_{i-1/2,j} = - c_D \alpha_D v_D.
\end{align*}
We will need to solve two of these transverse solvers for the Euler equations as shown in the grid in \Fig{trans}. 
We will only consider the up-going fluctuation of the transverse solver at $(i,j+1/2)$ and the down-going fluctuation of the
solver at $i,j-1/2$. This yields the full fluctuations as
 \begin{gather*}
\begin{gathered}
\B^+\A^+\Delta Q_{i-1/2,j} = 
 \frac{c_3 \left( c_2 \A^+_{1} + n_{3x} \A^+_{2} + n_{3y} \A^+_{3}  \right) }{c_3 + c_2}
\left[\begin{array}{c}  1 \\ n_{3x} c_3 \\ n_{3y} c_3  \\ 0 \end{array} \right], \\
\B^-\A^+\Delta Q_{i-1/2,j} = 
 \frac{-c_1\left( c_2 \A^+_{1} - n_{2x} \A^+_{2} - n_{2y} \A^+_{3}  \right)}{c_1 + c_2}
\left[\begin{array}{c}  1 \\ -n_{2x} c_1 \\  -n_{2y} c_1 \\ 0   \end{array} \right],
\end{gathered}
\end{gather*}
where $c_1,c_2$ and $c_3$ are the speeds of sound in cells $(i,j-1)$, $(i,j)$ and $(i,j+1)$ respectively, 
the normals $\hat{n}_3$ and $\hat{n}_2$ are the normals to the upper
edge and the lower edge, as shown in \Fig{trans}, and the 
non-acoustic fluctuations were neglected. The sound speeds are calculated with the pressure, density and the 
parameters of the Tammann EOS in the respective cell with $c=\sqrt{\gamma\frac{p+p_{\infty}}{\rho}}$. This 
is repeated analogously for the left going fluctuation $\A^-\Delta Q_{i-1/2,j}$ of the normal Riemann problem.
These transverse Riemann solvers were also implemented in the simulations shown in \Fig{maphighP}.

\subsection{Transmission based limiters}
\begin{figure}[t]
\centering
\includegraphics[width=0.25\paperwidth]{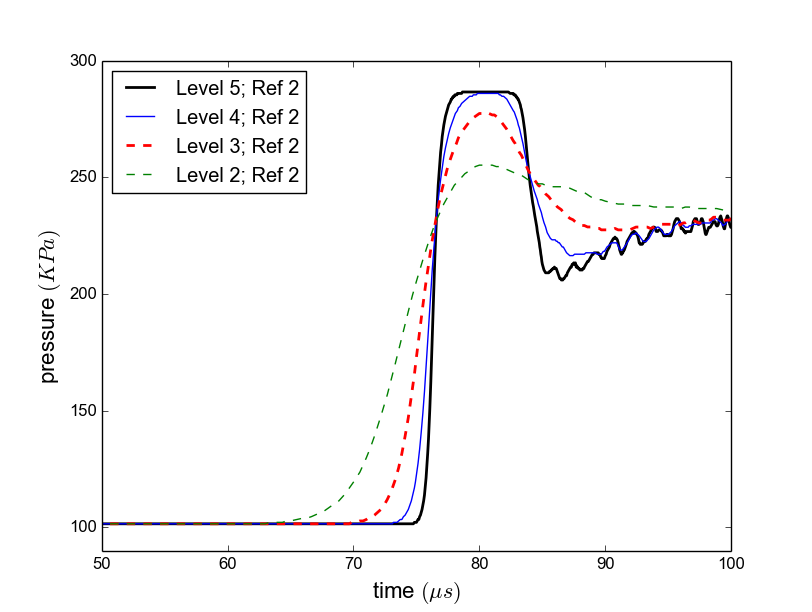} \ \ \
\includegraphics[width=0.25\paperwidth]{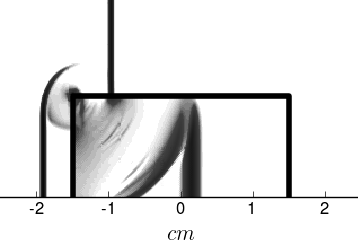} \\
\caption{In the first plot, we show a convergence study at a gauge at (-1cm,0). The curves shown are 
for four different AMR levels of refinement up to level 5, where each level doubles the resolution of 
the previous one. Oscillations are clearly seen in level 5 refinement. The second plot shows a schlieren plot for
the pressure where one can appreciate the oscillations produced at the corner of the interface.} 
\label{fig:trans-oscill}
\end{figure}

When the mesh is refined heavily by AMR, high-frequency unphysical oscillations appear in the water. 
Their wavelength scales with the mesh resolution, and they are hard to observe in the coarser grids 
due to numerical diffusion. These oscillations originate in the corner of the interface and they 
do not dissipate. This is caused by small errors produced by the Riemann solvers at the 
interface; these errors propagate in the normal and transverse direction. In the corner grid cell, these errors 
occur once when sweeping the solver on the grid horizontally 
and once again when sweeping vertically, producing oscillations. A sample of this phenomena can be observed in
\Fig{trans-oscill}, where we show the convergence study for a pressure gauge at (-1cm,0) and a schlieren
plot of the pressure that shows the oscillations being produced at the corner of the interface. The convergence
will be further studied in Section \ref{sec:Verif}.

\begin{figure}[t]
\centering
\includegraphics[width=0.25\paperwidth]{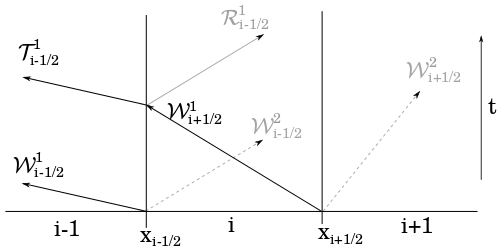} \ \ \
\includegraphics[width=0.25\paperwidth]{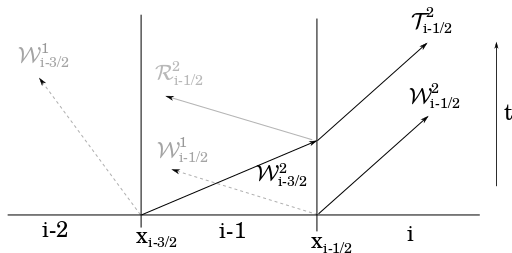} \\
\caption{Two diagrams are shown to illustrate the waves being compared in the different kind of limiters at the edge $i-1/2$, between
grid cells $i-1$ and $i$. The first diagram 
shows the waves that are involved in determining $\theta_{i-1/2}^1$ for the limiting behavior of $\W_{i-1/2}^1$. The second one
shows the waves involved in determining $\theta_{i-1/2}^2$ for the limiting behavior of $\W_{i-1/2}^2$. The notation is $\T$ for transmitted
waves and $\R$ for reflected ones.} 
\label{fig:translim}
\end{figure}

This issue can be improved by adjusting how the waves at the interface are limited. The limited waves 
from Eq. \ref{eq:Godunovpt2} are given by
\begin{align*}
 \widetilde{\W}_{i \pm 1/2}^p = \phi(\theta)_{i\pm 1/2}^p \W_{i \pm 1/2}^p,
\end{align*}
where $\phi(\theta)$ is the flux-limiter function \cite{randysrbook} and $\theta$ is a measurement of the
smoothness of the function. There are several ways to choose the $\theta$ parameter to limit the waves coming out of the edges at $i-1/2$. For 
a linear problem with two waves, where the 1-waves propagate to the left and the 2-waves to the right (like acoustics),
we can measure the smoothness $\theta$ by comparing the magnitude of adjacent waves. The corresponding $\theta$ parameters can be
obtained as $\theta_{i-1/2}^1 = \| \W_{i+1/2}^1 \|/\| \W_{i-1/2}^1 \|$ and $\theta_{i-1/2}^2 = \| \W_{i-3/2}^2 \|/\| \W_{i-1/2}^2 \|$,
see \cite{randysrbook}. In the case of nonlinear equations, the approach is similar; however, the eigenvectors of adjacent waves are 
no longer co-linear in phase space across adjacent cells, so we need to do a projection into the corresponding eigenvectors. 
For the nonlinear case, the $\theta$ parameter is given by
$\theta_{i-1/2}^1 = (\W_{i+1/2}^1 \cdot \W_{i-1/2}^1)/(\W_{i-1/2}^1 \cdot \W_{i-1/2}^1)$ and 
$\theta_{i-1/2}^2 = (\W_{i-3/2}^2 \cdot \W_{i-1/2}^2)/(\W_{i-1/2}^2 \cdot \W_{i-1/2}^2)$, see \cite{randysrbook}. 
The diagrams in \Fig{translim} give some visual intuition into which waves we are comparing. This is the standard
implementation in Clawpack \cite{clawpack}. 

In the case where there is a big jump in the parameters across an interface, the eigenvectors of a wave on different sides of 
the interface are significantly different. In this case, it is more appropriate to separate one of the adjacent waves into its transmitted
and reflected component, as if it actually had crossed the interface, and use the transmitted wave to limit the other adjacent wave. 
For instance assume the interface is at the edge $i-1/2$ shown in \Fig{translim}, the 
original limiter compares the projection of $W_{i+1/2}^1$ (into the corresponding eigenvector at $i-1$) 
with $W_{i-1/2}^1$ to limit $W_{i-1/2}^1$. However, 
if the interface has a big jump in the parameters, it is better to separate $W_{i+1/2}^1$ into its reflected and transmitted components and 
compare the transmitted component of the wave $T_{i-1/2}^1$ with $W_{i-1/2}^1$ to limit $W_{i-1/2}^1$. These type of limiters 
are called transmission based limiters, originally developed in \cite{fogarty1999high} for acoustics equations in heteregeneous media. 
In this case, the $\theta$ parameters are given by
\begin{align}
 \theta_{i-1/2}^1 = \frac{\| \T_{i-1/2}^1\|}{\| \W_{i-1/2}^1\|} \ \ \ \theta_{i-1/2}^2 = \frac{\| \T_{i-1/2}^2\|}{\| \W_{i-1/2}^2\|},
 \label{eq:translim}
\end{align}
where the transmitted waves $\T_{i-1/2}^{(1,2)}$ are as shown in \Fig{translim} correspondingly. This requires calculating the transmitted waves,
which might follow different procedures depending on the equations we are using.

In this section, we extend the methods in \cite{fogarty1999high} for acoustic equations to limit the acoustic waves in 
Euler equations. In order to do so, lets recall that we can rewrite the one-dimensional acoustic 
equations in terms of the density and the momentum \cite{delrazo2015_01,randysrbook},
 \begin{gather}
\label{eq:linacoust2}
\begin{gathered}
    \left[\begin{array}{c} \rho \\   \rho u \end{array} \right]_t +
    \left[\begin{array}{c c } 0 & 1  \\  
                        c^2 & 0 \end{array} \right]
    \left[\begin{array}{c} \rho \\ \rho u \end{array} \right]_x
= 0,
\end{gathered}
\end{gather}
where $c$ is the sound speed, the eigenvalues of the system at a cell interface are the left and right sound speeds, 
$\lambda_{1,2} = -c_{i-1}, c_i$, and the corresponding eigenvectors $r^1_{i-1} = [1,-c_{i-1}]$ and $r^2_i=[1,c_i]$. As we assume 
different materials accross the interface $c_{i-1}\neq c_i$. Following the first diagram of \Fig{translim} and Eqs. \ref{eq:translim}, in order 
to calculate $\theta_{i-1/2}^1$, we need to know $\T_{i-1/2}^1$, which is the transmitted wave from wave 
$\W_{i+1/2}^1$ coming from cell $i$ to cell $i-1$. In order to do so, we first write the wave $\W_{i+1/2}^1$ in 
terms of the corresponding eigenvector $\W_{i+1/2}^1=\alpha_{i+1/2}^1 r_i^1$, which we already know from solving the Riemann problem, see 
\cite{randysrbook}. Then we decompose it into the eigenvectors of the corresponding two cells to obtain the transmitted and reflected contributions,
\begin{align*}
 \alpha_{i+1/2}^1 \left[\begin{array}{c } 1  \\  
                  -c_i \end{array} \right] = 
                  \beta_{i+1/2}^1 \left[\begin{array}{c } 1  \\  
                  -c_{i-1} \end{array} \right] +
                  \beta_{i+1/2}^2 \left[\begin{array}{c } 1  \\  
                  c_i \end{array} \right].
\end{align*}
This yields two equations with two unknowns, so we can solve for the $\beta_{i+1/2}^1$,
\begin{align*}
 \beta_{i+1/2}^1 = \alpha_{i+1/2}^1\frac{2 c_i}{c_{i-1}+c_i}.
\end{align*}
This quantity multiplied by the eigenvector $r_{i-1}^1$ corresponds to the transmitted wave. With this information,
and using that $\W_{i-1/2}^1 = \alpha_{i-1/2}^1 r_{i-1}^1$, we can now calculate the $\theta$ parameter,
\begin{align*}
 \theta_{i-1/2}^1 = \frac{\| \T_{i-1/2}^1\|}{\| \W_{i-1/2}^1\|} = \frac{\alpha_{i+1/2}^1}{\alpha_{i-1/2}^1}\left(\frac{2 c_i}{c_{i-1}+c_i}\right), \\
 \theta_{i-1/2}^2 = \frac{\| \T_{i-1/2}^2\|}{\| \W_{i-1/2}^2\|} = \frac{\alpha_{i-3/2}^2}{\alpha_{i-1/2}^2}\left(\frac{2 c_{i-1}}{c_{i}+c_{i-1}}\right),
\end{align*}
where $\theta_{i-1/2}^2$ is calculated in the same manner by following the second diagram from \Fig{translim}. The sound speeds 
can be obtained from the Tammann EOS by using Eqs. \ref{eq:sos} and \ref{sndsp}. Also note the limiters work on the waves in the computational
domain, so it is not necessary to do any additional adjustments when using a mapped grid. 

These limiters greatly improve the observed oscillations as shown in the first plot of \Fig{AMRconv} where the level 5 refinement no longer
shows significant oscillations. Note these limiters are approximate since we are using the acoustic equations 
rewritten in terms of density and momentum to limit the Euler equations, and they don't fully suppress the oscillations in higher refinement levels 
as we will see in the next Section.

\section{Verification}
\label{sec:Verif}

A verification study for the one-dimensional case was performed in a 
previous work \cite{delrazo2015_01}. In that work, we verified that 
the finite volume methods coupled with the hybrid Riemann HLLC-exact 
Riemann solver for the Euler equations with a Tammann EOS converge to the 
correct solution for a simple model problem. However, the exact analytic solutions 
of Riemann problems for Euler equations are only available
in one dimension, so we restricted our verification to a one-dimensional 
test. Nonetheless, as the Riemann problem is still the key ingredient of higher-dimensional 
numerical methods,  the analysis from \cite{delrazo2015_01} is still relevant for the 
two-dimensional extension of the algorithm.

\begin{figure}[t]
\centering
\includegraphics[width=0.25\paperwidth]{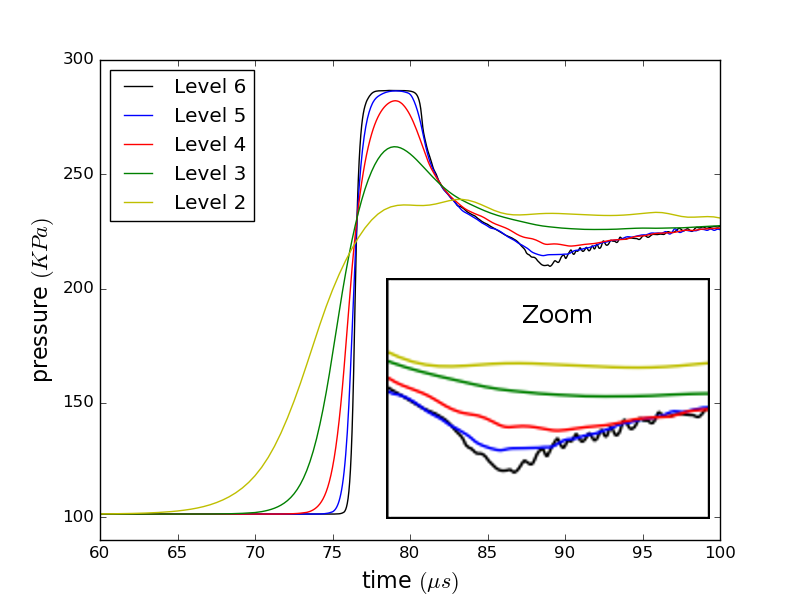}
\includegraphics[width=0.25\paperwidth]{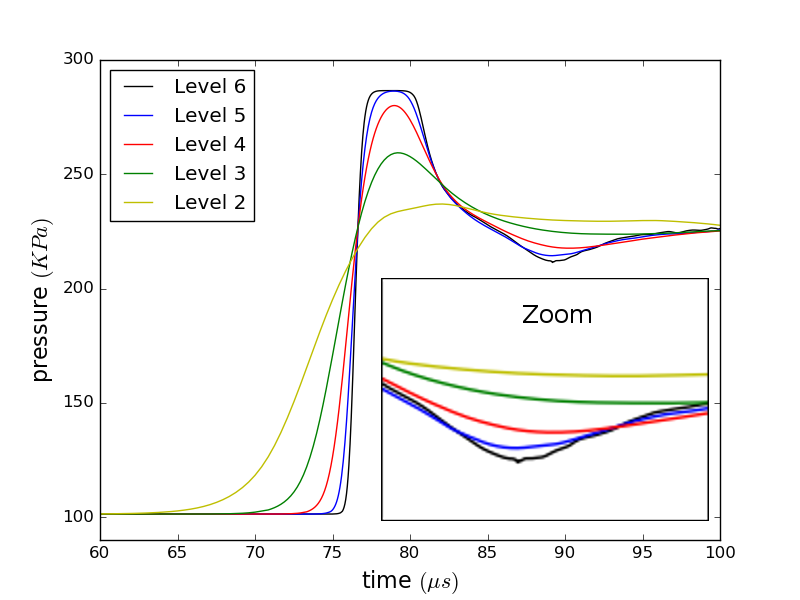} \\
\includegraphics[width=0.2\paperwidth]{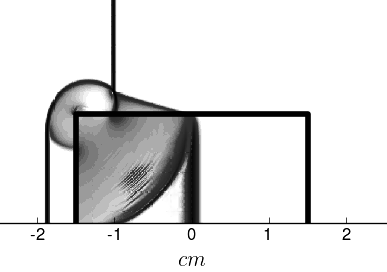} \ \ \
\includegraphics[width=0.25\paperwidth]{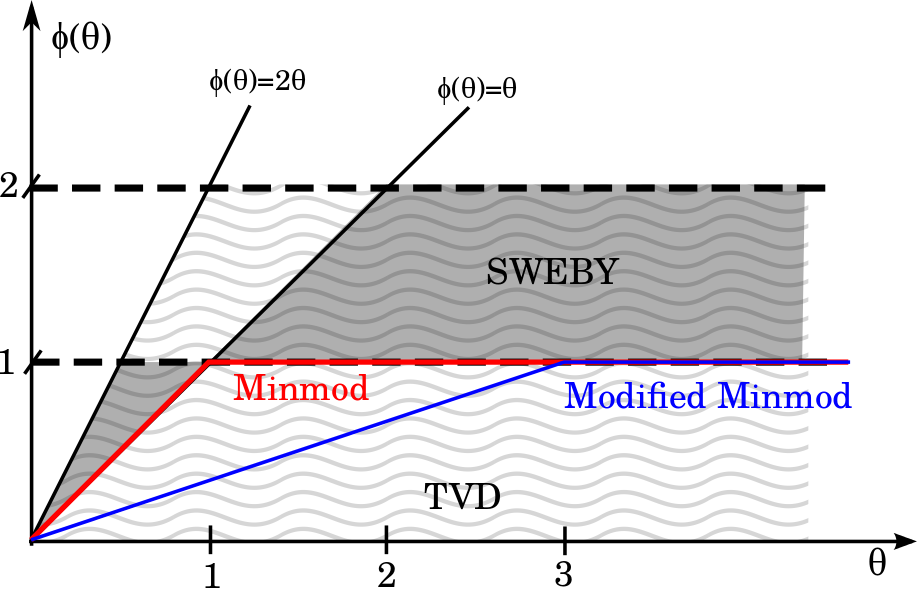}
\caption{The first two plots show the convergence tests at a gauge in (-1cm,0) for the two dimensional axisymmetric model 
with AMR on a Cartesian grid. The curves are shown for different levels of refinement allowed in AMR, where each level
doubles the resolution of the previous one. In the first plot one can appreciate numerical high frequency oscillations 
in the finer grids; however, this are almost fully supressed in the second figure by using a more diffusive limiter, the 
modified minmod limiter. The third plot shows a schlieren plot of the pressure with level 6 refinement when using the original
minmod limiter. It shows high-frequency oscillations propagating from the corner that do not dissipate. The fourth plot 
shows the TVD region (wavy lines) and the Sweby region (shaded) \cite{randysrbook} as well as the corresponding minmod and 
modified minmod limiter.} 
\label{fig:AMRconv}
\end{figure}

In addition to the verification study presented in \cite{delrazo2015_01}, in this work we will
provide a convergence test for the two-dimensional axisymmetric model. As there are 
no exact solutions for the two-dimensional equations, 
the convergence test only shows the numerical algorithm converges to a solution as the 
mesh is refined. The convergence tests were performed using several gauges for the Cartesian 
grid simulations of \Fig{cartesian}. In the first plot in \Fig{AMRconv}, we show the convergence test for 
the gauge at (-1cm,0). Note the appearance of high-frequency oscillations in the most refined level (level 6) 
even after applying the transmission-based limiters. The plot at the bottom-left
of \Fig{AMRconv} shows these oscillations for the finest grid in a schlieren pressure plot. 
The origin of this oscillations is the same as before. 

These oscillations can be suppressed by adding some numerical viscosity to the water material. This is not
entirely unphysical since the water is a viscous media. In order to do so, we implement a new limiter for the water grid cells, which we 
refer to as modified minmod. The original minmod limiter uses the flux-limiter function 
$\phi(\theta) = \mathrm{minmod}(1,\theta)$ \cite{randysrbook}. The minmod limiter is the most dissipative second-order
total variation diminishing (TVD) limiter. This is shown in the flux-limiter function plot at the bottom-right 
of \Fig{AMRconv}. The region covered in wavy lines is the region where the limiter can be TVD, and the shaded region shows
the Sweby region where limiter can be second-order accurate; the corresponding flux-limiter function for the minmod limiter 
is shown too. In order to add more numerical viscosity, we use a modified minmod limiter 
$\phi(\theta) = \mathrm{minmod}(1,\theta/3)$. Although we lose second order accuracy for the Euler Equations, this limiter
still provides physical solutions due to water viscosity. The scaling factor within the flux-limiter function ($1/3$) was 
chosen to be as close to 1 as possible to keep as much overlap with the Sweby region as possible while also supressing the oscillations; 
this parameter can be easily adjusted in the code, which is available in \cite{TBI-zenodo2}. The resulting convergence study 
after applying the modified minmod limiter can be appreciated in the second plot of \Fig{AMRconv}, where we can observe 
the oscillations were suppressed and that our method converges. Analogous results were obtained for the other gauges.

\section{Discussion}
\label{sec:TBIdisc}

We developed a two-dimensional axisymmetric shock-capturing high-resolution numerical model to study shock wave 
dynamics when crossing a fixed interface between a compressible fluid (air) and an almost incompressible material (water). 
These methods have been designed to complement TBI and other biomedical experiments performed in a shock tube. The common setup in these
experiments consists of a shock wave traveling through air and impacting a plastic container. The container is usually very 
thin, and it is often filled with an aqueous solution where 
the biological sample is placed. In our computational simulations, the container is modeled as an interface fixed in space. 
The aim of these methods and simulations is to provide experimentalists measurements of relevant variables inside the container, 
like pressure, that would otherwise be very difficult to obtain experimentally. This can help us understand better the
on-going physical dynamics that experimental samples in specifc experiments undergo and explain possible damage mechanisms. 
It should be noted the methods developed here can be extended to other scenarios.

We first provided the one-dimensional methods employed in detail and their implementation into Clawpack \cite{clawpack}.
In Appendix \ref{sec:TBI:compexp} we show that there is not a significant difference between the transmitted 
shock wave when removing the thin plastic interface separating air and water. Furthermore, we observed an 
amplification and elongation of the shock wave. This effect is accounted for by the different material compressibility. 
The amplitude of the initial pressure wave in the air increased in $54 \%$ when measured in the water. This amplification 
effect was highly relevant in the injury mechanisms studied in \cite{delrazo2015_01}, and it generally occurs 
when passing from air to water or a solid material. 

The methods were extended to two dimensions and implemented on a mapped grid, which allows more complicated 
interface geometries as long as the mapping is provided. We provided as a proof of concept a circular inclusion 
mapping, which maps the rectangular interface into a circular one. In the axisymmetric case, this mapping models 
a spherical interface. In addition, the algorithms 
were adapted to work with AMR to increase resolution and efficiency of the code. Additional mathematical work has
to be performed to improve the accuracy and stability of the numerical method. Transverse Riemann solvers for the 
mapped grid were developed to improve the accuracy. Transmission-based limiters and the minmod modified limiter 
were implemented at the interface and in the water to suppress numerical oscillations at heavily refined AMR patches.
A more primitive version of these methods was already employed in a specific mild TBI 
application \cite{delrazo2015_01}, and we expect they can be extended and used in new applications. 
The three-dimensional model with a spherical shell inclusion could be specifically useful in TBI applications 
to model an idealized skull of a mouse inside the shock tube or a human head exposed to a shock, a problem of much
interest to the TBI community as shown by some previous studies \cite{
anderson2003modern,ho2007dynamic,kleiven2007predictors,moss2009skull,ruan1993finite,shreiber1997vivo,
sundaramurthy2012blast,takhounts2003development,takhounts2008investigation,
taylor2009simulation,wakeland2005computer,zhou1997viscoelastic} among others. 
The code where all these methods are implemented is available with a BSD license \cite{TBI-zenodo2}.

\appendix
\numberwithin{equation}{section}

\section{One dimensional computational experiments}
\label{sec:TBI:compexp}
In this section, we simulate the one-dimensional Euler equations \ref{eq:Eulercyltamman} with the numerical methods from Section \ref{sec:nummethods}. We explore the question of whether
a thin plastic interface separating gas and liquid in a shock tube experiment can be ignored
in computational experiments, specifically whether the magnitude of the shock wave transmitted from the gas 
to the liquid is insensitive to the intervening layer of plastic.  
In the laboratory experiments that motivated this work, the walls of the plastic transwell container
are thin relative to the dimensions of the interior, and the computations presented in 
\cite{delrazo2015_01} were simplified by omitting the plastic layer entirely.  
Here we justify that approximation by considering a simple
one-dimensional model of a shock wave passing through layers of air-plastic-water. 
This will provide insight on the behavior of the shock wave; it will show what parameters are the most relevant, 
and it will show that it is not necessary to include the thin plastic interface in the computational model. We will 
begin with a one-dimensional air-plastic-water interface using Euler equations, 
and then we compare these results to the simpler one-dimensional air-water interface, omitting the plastic layer. 
We will further verify our results with an analytic calculation for a thin interface in linear acoustics. A preliminary 
version of this study can also be found in the conference proceedings \cite{delrazo2014}. 

\subsection{Air-plastic-water interface}
\label{sec:1Dapwinterface}
We begin by studying an air-plastic-water interface with Euler equations \ref{eq:Eulercyl} in one dimension,
\begin{gather}
\label{eq:Eulercyl1D}
\begin{gathered}
  \fr{\p}{\p t}
    \left[\begin{array}{c} \rho \\   \rho u \\    E   \end{array} \right] 
+ \fr{\p}{\p x} 
    \left[\begin{array}{c} \rho u \\   \rho u^2 + p \\   u(E+p)   \end{array} \right] 
=  0.
\end{gathered}
\end{gather}

\begin{figure}[H]
\centering
\includegraphics[width=0.38\paperwidth]{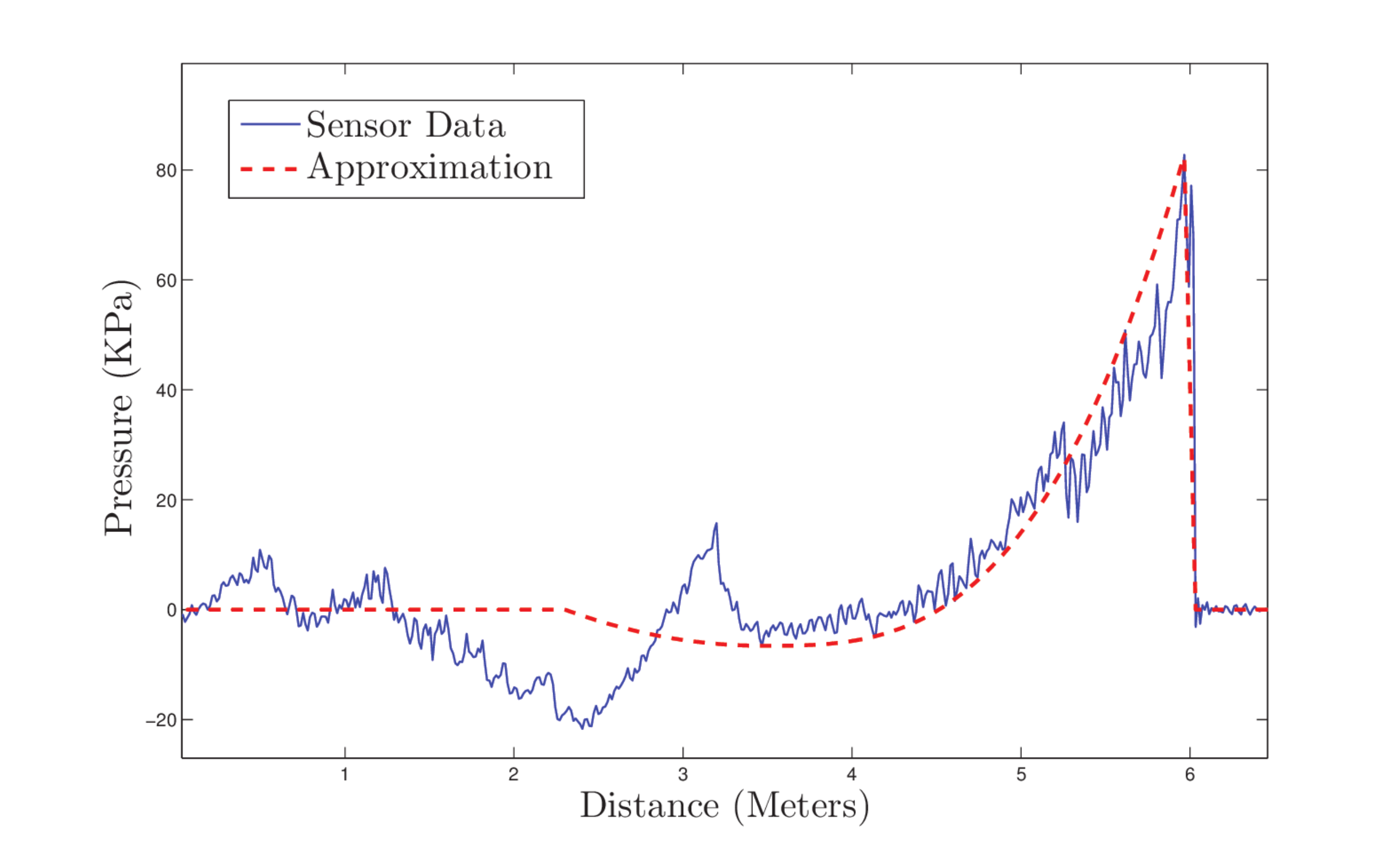}
\caption{The shock wave form obtained from a sensor inside a shock tube
is shown as the solid thin line. The coarse approximation 
to be used as an initial condition in our simulation is shown with a dashed line. An 
average speed of sound of $c=344$ $m$/$s$ is assumed. These figure was obtained from 
\cite{delrazo2015_01}.}
\label{fig:dataTBI}
\end{figure}

\begin{figure*}[t]
  \centering
  \subfigure[]{ \includegraphics[width=0.3\textwidth]{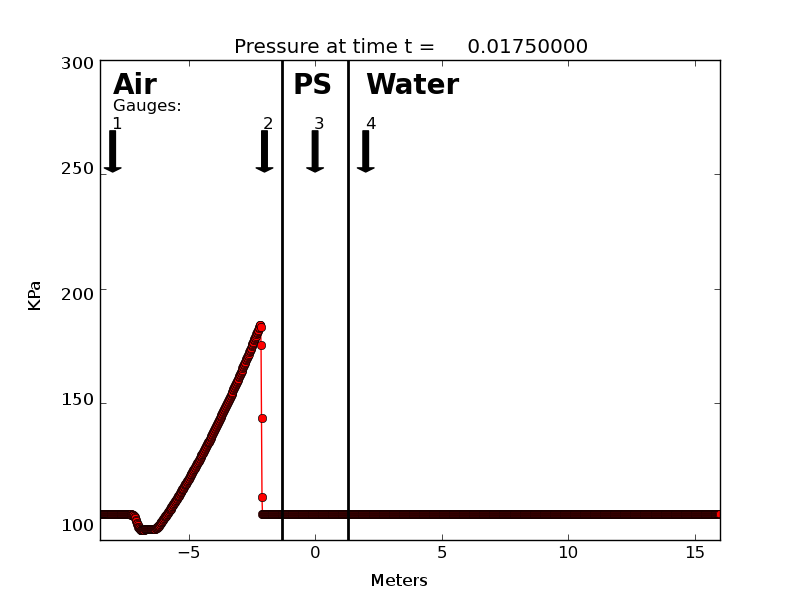} 
      \label{fig:TBIanim1}}
  \subfigure[]{ \includegraphics[width=0.3\textwidth]{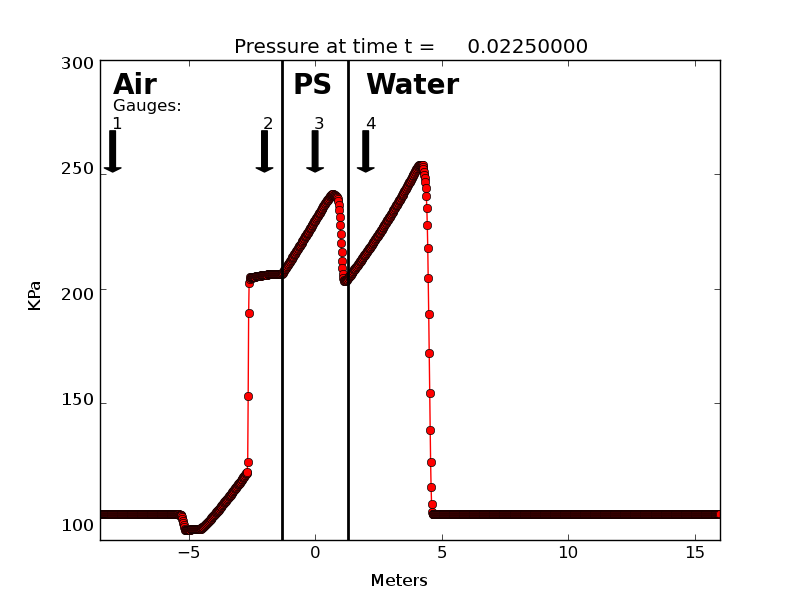} 
      \label{fig:TBIanim2} }
  \subfigure[]{ \includegraphics[width=0.3\textwidth]{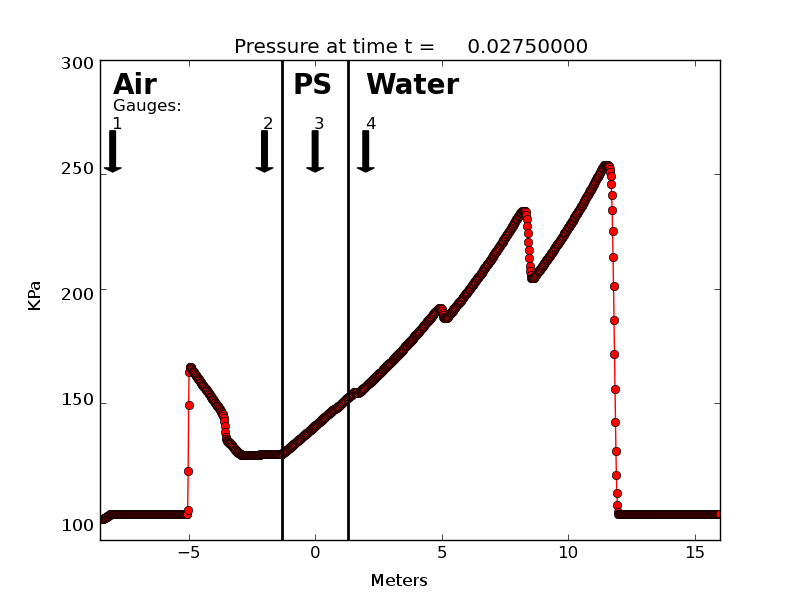} 
      \label{fig:TBIanim3} }
  \caption{Shock wave crossing the air-plastic-water interface at different times. The arrows indicate the position
of the 4 gauges that measure the pressure as a function of time. The gauges are numbered from left to right, 
and the plastic interface width for this case is $2.6$m }
  \label{fig:TBIanim}
\end{figure*}

The first step is to input the right initial conditions into our simulation. The actual form of the initial shock wave 
traveling through the shock tube was obtained experimentally; the amplitude can be varied
in the shock tube and in our computational simulations. The
sensor outputs pressure amplitude as a function of time. Assuming an average
speed of sound in air, it can be converted to a function of distance as shown in
Figure \ref{fig:dataTBI}. The shape can be broadly approximated by an idealized
shock wave (dashed line in Figure \ref{fig:dataTBI}). This approximated shape of
the shock wave is introduced as the initial condition in the simulation, where a scaling factor
is used to scale the amplitude; however,
this is not a trivial procedure since we must input the density, momentum, and energy, and we only have 
the pressure. Using the isentropic EOS, the ideal gas EOS and the expression for the speed of sound, an 
educated guess for the initial condition in terms of the pressure is given far away from the transwell. 
This initial condition is then modified until we obtain the desired amplitude and shape of the shock wave 
front. The resulting shape of the shock wave before hitting the interface can be seen in 
Figure \ref{fig:TBIanim1}, where the scaling factor, in this case, was chosen arbitrarily. The pressure is 
measured in KPa with an ambient base pressure of $1\text{ATM}=101.325\text{KPa}$. The same procedure was used for the 
two-dimensional simulations.

The one-dimensional equations with the pair of interfaces are solved using the methods 
mentioned in Section \ref{sec:nummethods}. The different materials are modeled using 
different parameters for the Tammann EOS, see Section \ref{sec:EOS}. The choice of 
parameters is shown in Table \ref{tab:param}.

The solution of the shock wave crossing the interfaces between air, plastic, and water 
at different times is shown in Figure \ref{fig:TBIanim}. It can be observed that every time the shock
wave hits an interface, part of the wave is reflected and part of it is transmitted. This effect can occur 
multiple times depending on how the interfaces are set up. We can also observe that 
the amplitude of the shock wave increases as it passes from air to plastic and decreases 
when passing from plastic to water. This effect is due to the
continuity of pressure and the change in compressibility. In order to keep the pressure at the interface 
continuous, the transmitted wave amplitude has to be the same as the sum of the incident wave and 
the reflected wave. When the compressibility is very high in the adjacent material, the interface 
will behave similarly to a solid wall. In this case, since the reflected wave will have an amplitude 
almost equal to the incident wave, the transmitted wave could have an amplitude almost twice 
as big as that of the incident wave. This explains why the pressure jump can increase or decrease when
crossing an interface. Even for the one-dimensional case, we observe complex behavior due to interaction at 
the interface. These numerical simulations provide 
accurate insight in situations where simple intuition might be insufficient.  

In Figure \ref{fig:dataTBI}, we show from experimental data the initial shock wave profile in the air before hitting 
any interface; however, we are interested in the shape and amplitude of the shock wave in the water. 
In order to do so, we first need to know how important the plastic interface 
is in our model. Computationally, the plastic interface is hard to model because the 
width of the plastic is very small ($~mm$) in comparison 
to the characteristic length of the experiment (length of the transwell \cite{delrazo2015_01}). 
The following experiment explores how the width of the plastic 
interface affects the shock wave profile. Additionally, we show an accurate model can be obtained even 
when completely ignoring the plastic interface.

\begin{figure*}
  \centering
  \subfigure[]{ \includegraphics[width=0.3\textwidth]{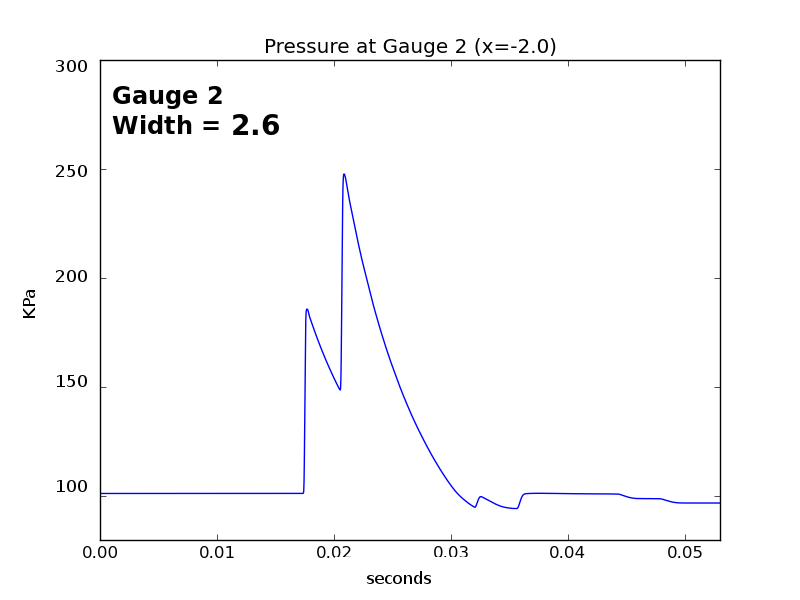} 
      \label{fig:TBIgaugea} }
  \subfigure[]{ \includegraphics[width=0.3\textwidth]{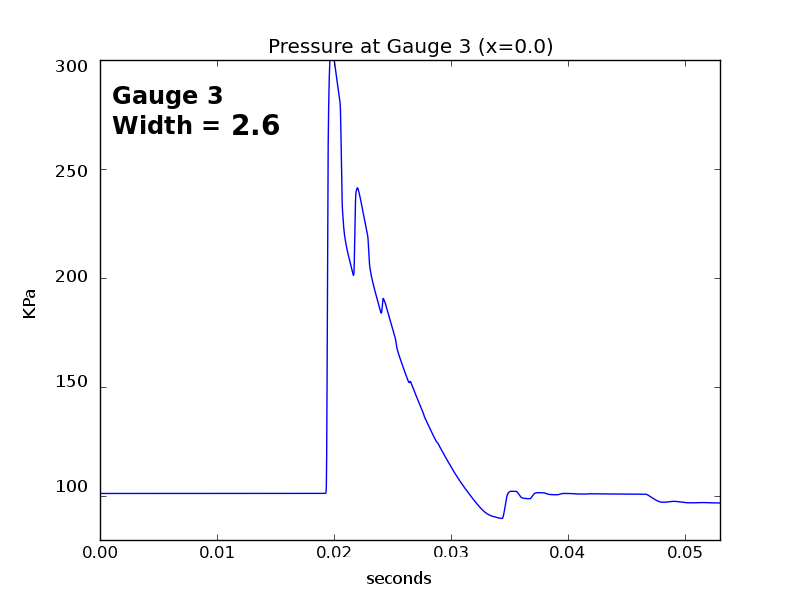} 
      \label{fig:TBIgaugeb} }
  \subfigure[]{ \includegraphics[width=0.3\textwidth]{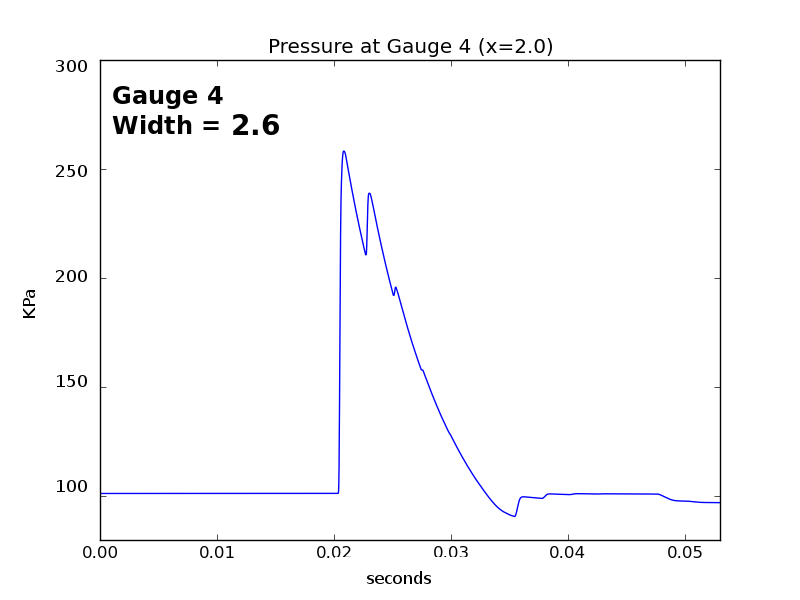} 
      \label{fig:TBIgaugec} }
  \subfigure[]{ \includegraphics[width=0.3\textwidth]{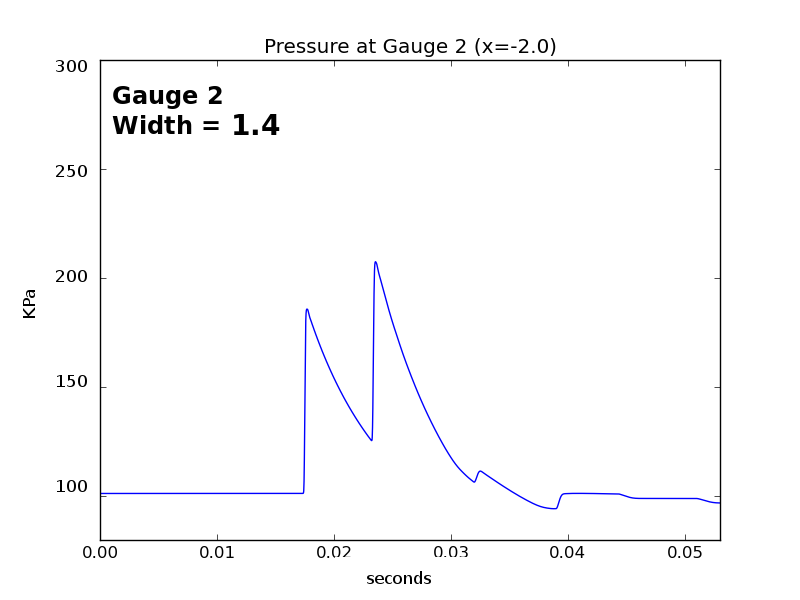} 
      \label{fig:TBIgauged} }
  \subfigure[]{ \includegraphics[width=0.3\textwidth]{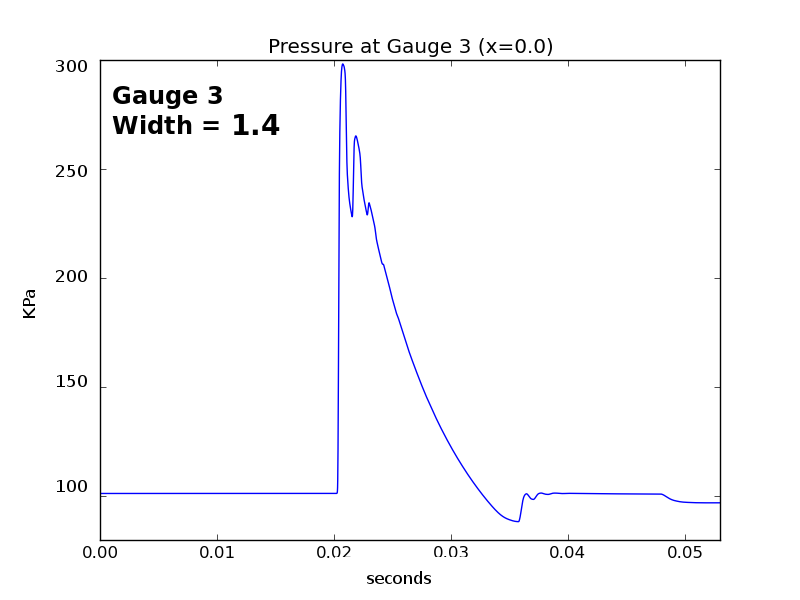} 
      \label{fig:TBIgaugee} }
  \subfigure[]{ \includegraphics[width=0.3\textwidth]{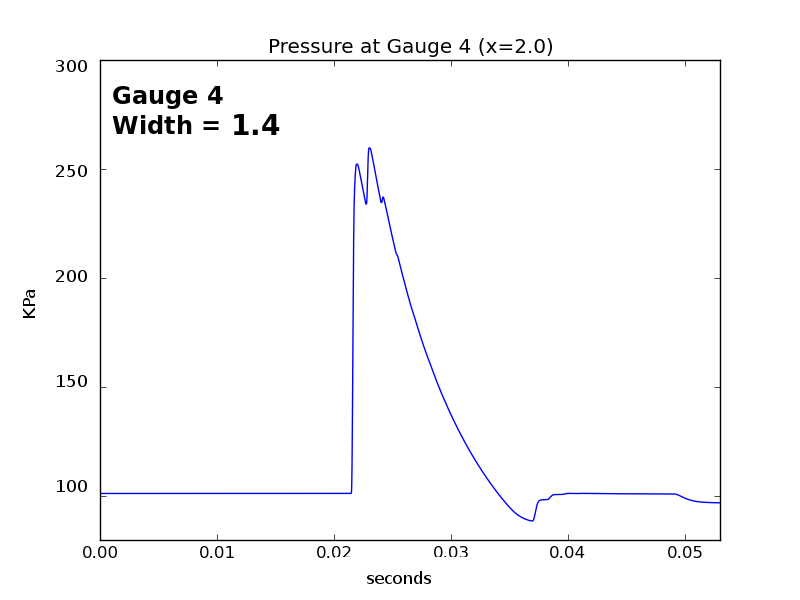} 
      \label{fig:TBIgaugef} }
  \subfigure[]{ \includegraphics[width=0.3\textwidth]{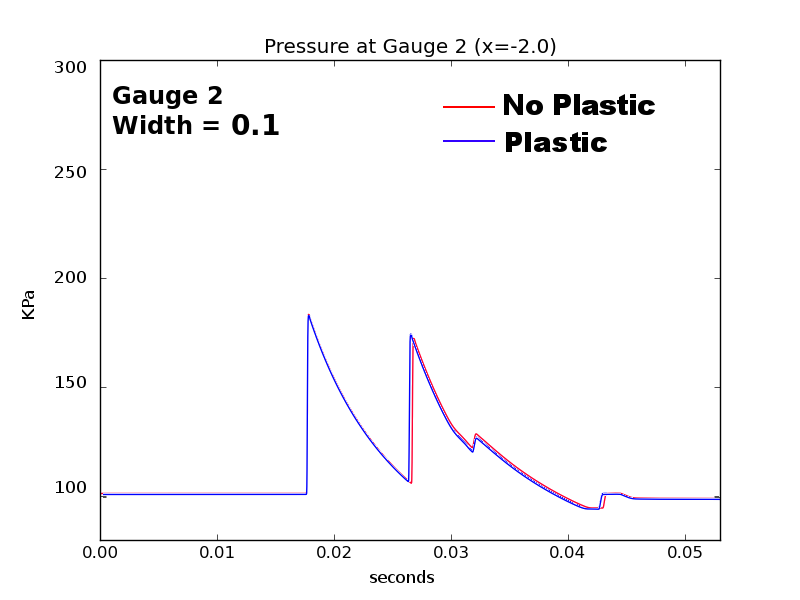} 
      \label{fig:TBIgaugeg} }
  \subfigure[]{ \includegraphics[width=0.3\textwidth]{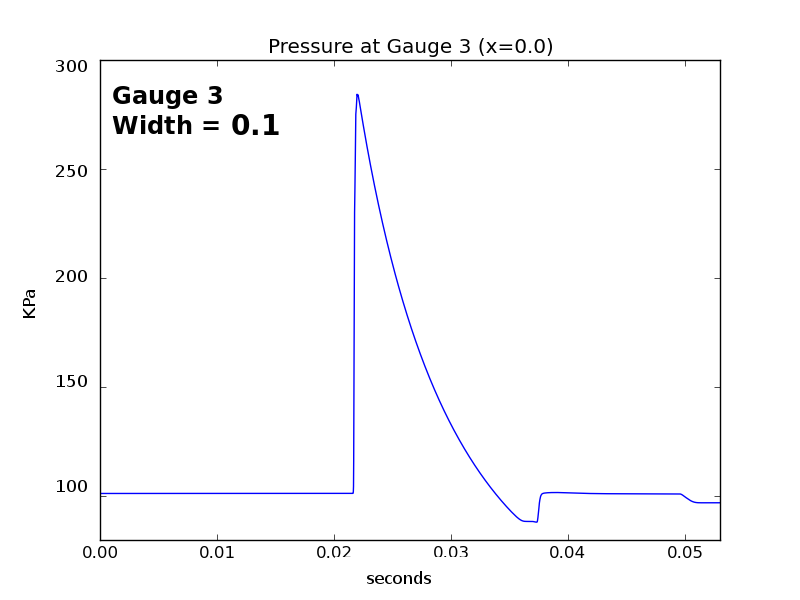} 
      \label{fig:TBIgaugeh} }
  \subfigure[]{ \includegraphics[width=0.3\textwidth]{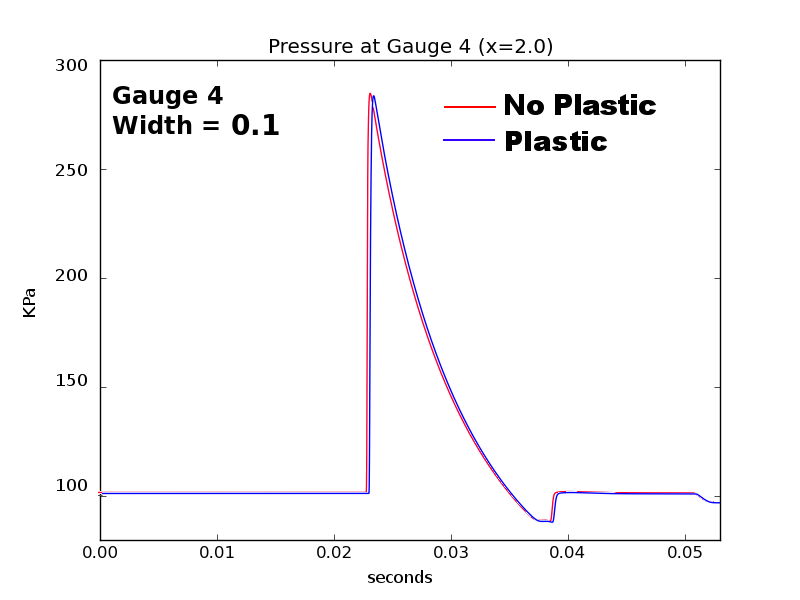} 
      \label{fig:TBIgaugei} }
  \caption{Pressure (KPa) gauge plots as a function of time (seconds). Each row of figures shows the 
  three gauge plots for three different widths ($2.6$m, $1.4$m and $0.1$m) of the plastic interface, 
  as shown in Table \ref{tab:widthamp}. The plots (g) and (i) for gauge 2 and 4 also show the pressure gauge 
  plots when there is no plastic interface at all; the difference is almost unnoticeable. Also note the red line in
  Figure \ref{fig:TBIgaugeg} is completely overlapped by the blue line before the reflected shock appears; this is 
  because the solutions between thin plastic and no plastic are exactly the same before interacting with the interface.}
  \label{fig:TBIgauges}
\end{figure*}

The maximum amplitude of the pressure profile was measured at gauges $2,3$ and $4$ of 
Figure \ref{fig:TBIanim} for different widths of the plastic
interface. The plastic is always assumed to be centered at $x=0$. The results are presented 
in Table \ref{tab:widthamp}. In Figure \ref{fig:TBIgauges}, the full pressure profiles as a 
function of time are shown at the three gauges for three of the plastic widths shown in Table \ref{tab:widthamp}. 

\begin{table*}
  \centering
  \begin{tabular}{ l || c  c  c  c }
    \hline
     Width (m) & Initial(KPa)  & Gauge 2 (KPa) & Gauge 3 (KPa) & Gauge 4 (KPa) \\ \hline
    2.6        & 184.06    & 247.76    & 305.88    & 258.24    \\ 
    1.4        & 184.06    & 207.53    & 298.19    & 259.71    \\ 
    0.6        & 184.06    & 187.72    & 283.72    & 274.90    \\ 
    0.2        & 184.06    & 183.31    & 282.34    & 280.18    \\   
    0.1        & 184.06    & 183.31    & 284.29    & 283.55    \\
    0.0        & 184.06    & 184.40    & -         & 284.26
  \end{tabular}
  \caption{The maximum amplitude measured at three pressure gauges for different 
  widths of the plastic interface. The initial shock wave is the same for all cases, and the 
  gauge plots are placed before, inside and after the plastic interface as shown in Figure \ref{fig:TBIanim}. The last row corresponds to the air-water interface.}
  \label{tab:widthamp}
\end{table*}

The results in Table \ref{tab:widthamp} and Figure \ref{fig:TBIgauges} show the maximum amplitude at gauge 2 is reduced 
as the plastic width is decreased. Not surprisingly, this is a consequence of having less interference with the 
reflected shock wave, since the gauge is farther away from the interface as the plastic width is reduced. This effect 
is clearly shown in Figures \ref{fig:TBIgauges}a, \ref{fig:TBIgauges}d, \ref{fig:TBIgauges}g. The maximum amplitude at 
gauge 3 is somewhat diminished at first; however, it seems to be reaching a plateau around $280.0 \text{KPa}$. The 
behavior at gauge 3 is not trivial; the shock wave bounces back and forth several times, interfering with itself 
constantly. In Figures \ref{fig:TBIgauges}b, \ref{fig:TBIgauges}e, \ref{fig:TBIgauges}h, we can see the interference 
becomes so fast that the pressure profile in the plastic seems to converge to a shock wave shape as the plastic width 
is reduced. At gauge 4, we can observe the interference 
between the set of transmitted shock waves generated by the back and forth reflections within the plastic interface. 
As the plastic width is reduced, the time elapsed between the transmitted shock waves is reduced and the interference 
increased. Nonetheless, when the plastic width is very small, the interference becomes so fast that the pressure 
profile seems to converge again to a shock wave shape, as shown in Figures \ref{fig:TBIgauges}c, 
\ref{fig:TBIgauges}f, \ref{fig:TBIgauges}i. Furthermore, note the difference in the shock wave shape in 
Figures \ref{fig:TBIgauges}h, \ref{fig:TBIgauges}i is almost unnoticeable. 
It almost seems like the shock wave is only crossing one interface instead of two. This motivates the 
next experiment.

\subsection{Air-water interface}
\begin{figure}[H]
  \centering
  \subfigure[]{ \includegraphics[width=0.3\textwidth]{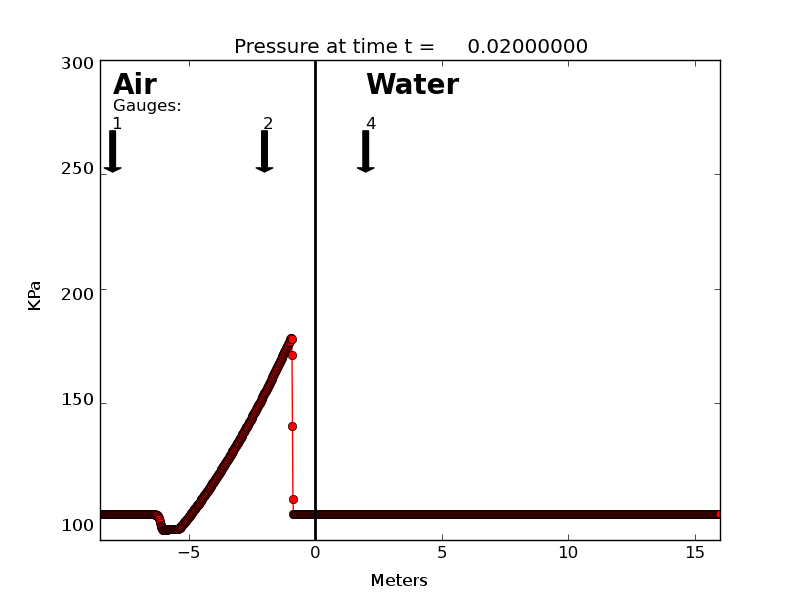} 
      \label{fig:TBIanimB1}}
  \subfigure[]{ \includegraphics[width=0.3\textwidth]{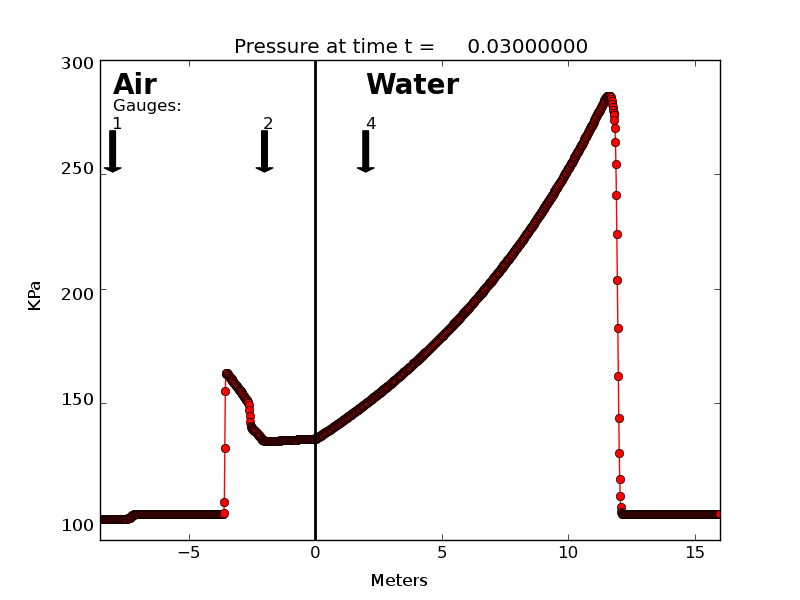} 
      \label{fig:TBIanimB2} }
  \caption{Shock wave before and after crossing the air-water interface. The arrows indicate the position
of the gauges that measure the pressure as a function of time. The gauges are the same as in Figure 
\ref{fig:TBIanim}. Gauge 3 was removed since  there is no plastic layer in this case.}
  \label{fig:TBIanimB}
\end{figure}

In reality, the plastic is so thin that is really unnoticeable on larger scales. 
Furthermore, as the plastic is almost an incompressible medium, one should expect 
it would transfer the shock wave infinitely fast without energy loss. Therefore, 
instead of the triple material interface, now consider only an air-water interface. 
The result of this simulation is shown in Figure \ref{fig:TBIanimB}. The gauge plots for 
gauge 2 and 4 are shown in Figures \ref{fig:TBIgauges}g and \ref{fig:TBIgauges}i, along with 
the thin plastic results. The maximum amplitude in each of these gauges is presented in the last row of Table \ref{tab:widthamp}.

Comparing the air-water interface results against the ones for the 
smallest plastic width in the air-plastic-water interface case, 
we can observe the percentage error in the maximum pressure amplitude of 
gauge 4 is of $0.38\%$. This is also obvious from the thin plastic and no plastic comparison
in Figures (\ref{fig:TBIgauges}g, \ref{fig:TBIgauges}i). 
This result allowed us to simplify higher dimensional air-plastic-water interface 
problem to a simpler air-water interface in the work \cite{delrazo2015_01}. Nonetheless, the presence of
the plastic is still modeled, since we force our interfaces to be fixed in space, just as a plastic container would 
force water to remain inside the container. 

\subsection{Air-Plastic-Water interface with linear acoustics}
\label{sec:exactac}
In order to further justify dropping the plastic layer, we consider the same situation 
of a thin intermediate layer for the case of linear acoustics.
In this case, we can compute the exact solution of the transmitted pressure
through the air-plastic-water interface as a function of the acoustic impedance of each material and the plastic width.
This can be derived from the fact that an acoustic wave with incident pressure jump $p_0$ on the left of an interface 
between medium $A$ (left) and $B$(right) produces a reflected and a transmitted wave with pressure jumps given by,
\begin{align*}
  p_T=p_0\frac{2Z_B}{Z_A + Z_B} \hspace{10mm} p_R=p_0\frac{Z_B-Z_A}{Z_A+Z_B},
\end{align*}
where $Z_k$ denotes the acoustic impedance of medium $k$. These relations can be easily derived from linear acoustics \cite{randysrbook}.
Now consider a one-dimensional air-plastic-water interface. With this
setup, there will be an infinite number of reflections in the plastic layer. 
The $N^{th}$ wave contribution to the transmitted wave in water is given by
\begin{align*}
  p_T^{N} = \frac{2Z_w}{Z_w+Z_p}\left(\frac{Z_a-Z_p}{Z_a+Z_p}\right)^{N-1}
                                     \left(\frac{Z_w-Z_p}{Z_w+Z_p}\right)^{N-1}
                 \frac{2Z_p}{Z_p+Z_a} p_0,
\end{align*}
where $Z_a$, $Z_p$, and $Z_w$ are the air, plastic and water impedances. 
Each transmitted wave increases the pressure behind the initial transmitted wave
slightly and the asymptotic final amplitude of the transmitted
wave is given by the sum of all these contributions,
\begin{gather*}
  p_{T}^{total} = \sum_{N=1}^{\infty} p_T^{N} = \\
\frac{4Z_w Z_p p_0}{(Z_w+Z_p)(Z_p+z_a)}\sum_{N=0}^{\infty} 
   \left(\frac{(Z_a-Z_p)(Z_w-Z_p)}{(Z_a+Z_p)(Z_w+Z_p)}\right)^{N}.
\end{gather*}
Summing this geometric series yields
\begin{align*}
  p_{T}^{total} = p_0\frac{2Z_w}{Z_w + Z_a}.
\end{align*}
When the plastic layer is very thin, this asymptotic value is quickly reached,
and we note that it is
exactly the same as if the plastic interface didn't exist. The transmission coefficient is 
the one computed directly from air into water. 
Note we assumed the pressure profile on the left was a constant $p_0$. However, this can be more complicated. It can have a 
decaying tail, in which case there will be interference from the tail in the reflected and transmitted waves. Nonetheless, 
assuming the plastic width is $w_0$, the time elapsed between two transmitted waves in the water interface is given by 
$\tau=2w_0/c_p$, where $c_p$ is the speed of sound in plastic. Therefore, as $w_0\rightarrow 0$, the elapsed 
time $\tau\rightarrow 0$. As a consequence, the interference from
the tail will also disappear and the plastic interface can be neglected without losing accuracy.

This calculation is an analytic result that shows that if the plastic interface is very thin in comparison to 
the experiment's characteristic length scales, the plastic interface can be neglected without losing much accuracy.

\bibliographystyle{siam}
\bibliography{extrarefs.bib,library.bib,biorefs.bib}

\end{document}